\documentclass{siamltex}

\usepackage{amsfonts}
\usepackage{epsfig}
\usepackage{amssymb}
\usepackage{amsmath}
\usepackage{algorithm}
\usepackage{algorithmic}
\usepackage{subfig}
\usepackage{graphicx}
\usepackage{color}
\usepackage[section] {placeins}
\usepackage{epstopdf}

\newcounter{ass_counter}

\newtheorem{assumption}[ass_counter]{Assumption}

\def\rc{\color{black}}

\def\hx{\hat{x}}
\def\P{\mathcal{P}}

\def\E{\mathbb{E}}
\def\PS{\mathcal{P}_{S}}
\def\PO{\mathcal{P}_{{\gamma \over \Lmax}g}}

\newcommand{\AsySPCD}{\mbox{\sc AsySPCD}}
\newcommand{\AsySCD}{\mbox{\sc AsySCD}}
\newcommand{\Lres}{L_{\mbox{\rm\scriptsize res}}}
\newcommand{\Lmax}{L_{\mbox{\rm\scriptsize max}}}
\newcommand{\beq}{\begin{equation}}
\newcommand{\eeq}{\end{equation}}
\def\eqnok#1{(\ref{#1})}
\newcommand{\R}{\mathbb{R}}

\def\hogwild{\mbox{\sc Hogwild!}}

\def\sjwcommentsolved#1{}
\def\jlcommentsolved#1{}

\def\citep#1{\cite{#1}}
\def\citet#1{\cite{#1}}

\title{Asynchronous Stochastic Coordinate Descent: Parallelism and Convergence Properties}

\author{Ji Liu\thanks{Department of Computer Sciences, University of Wisconsin-Madison, 1210 W. Dayton St., Madison, WI 53706-1685, US ({\tt ji.liu.uwisc@gmail.edu}). This author was supported in part by NSF Awards DMS-0914524 and DMS-1216318 and
ONR Award N00014-13-1-0129.}
\and Stephen J. Wright\thanks{Department of Computer Sciences, University of Wisconsin-Madison, 1210 W. Dayton St., Madison, WI 53706-1685, US ({\tt swright@cs.wisc.edu}). This author was supported in part by NSF Awards DMS-0914524, DMS-1216318, and IIS-1447449, ONR Award N00014-13-1-0129, AFOSR Award FA9550-13-1-0138, and Subcontract 3F-30222 from Argonne National Laboratory.}
}

\begin{document}

\maketitle

\begin{abstract}
We describe an asynchronous parallel stochastic proximal coordinate
descent algorithm for minimizing a composite objective function, which
consists of a smooth convex function added to a separable convex function.
In contrast to previous analyses, our model of asynchronous
computation accounts for the fact that components of the unknown
vector may be written by some cores simultaneously with being read by
others.  Despite the complications arising from this possibility, the
method achieves a linear convergence rate on functions that satisfy an
optimal strong convexity property and a sublinear rate ($1/k$) on
general convex functions. Near-linear speedup on a multicore system
can be expected if the number of processors is $O(n^{1/4})$. We
describe results from implementation on ten cores of a multicore
processor.
\end{abstract}



\pagestyle{myheadings}
\thispagestyle{plain}

\section{Introduction}
\label{sec:intro}

We consider the convex optimization problem
\begin{equation}
\label{eqn_mainproblem}
\min_{x} \, \quad F(x):=f(x) + g(x),
\end{equation}
where $f: \mathbb{R}^n \mapsto \mathbb{R}$ is a smooth convex function
and $g: \mathbb{R}^n \mapsto \mathbb{R}\cup \{\infty\}$ is a
separable, closed, convex, and extended real-valued
function. ``Separable'' means that $g(x)$ can be expressed as $g(x) =
\sum_{i=1}^n g_i((x)_i)$, where $(x)_i$ denotes the $i$th element of
$x$ and each $g_i: \mathbb{R}\mapsto \mathbb{R}\cup\{\infty\}$,
$i=1,2,\dotsc,n$ is a closed, convex, and extended real-valued
function.

Formulations of the type \eqref{eqn_mainproblem} arise in many data
analysis and machine learning problems, for example, the linear primal
or nonlinear dual formulation of support vector machines
\jlcommentsolved{Shall we mention that the bias term is removed? SJW: I
  think it is OK to omit this detail.}\citep{CortesVapnik95}, the
LASSO approach to regularized least squares, and regularized logistic
regression. Algorithms based on gradient and approximate / partial
gradient information have proved effective in these settings. We
mention in particular gradient projection and its accelerated variants
\citep{nesterov2004introductory}, proximal gradient \citep{WriNF08a}
and accelerated proximal gradient \citep{BeckT09} methods for
regularized objectives, and stochastic gradient methods
\citep{Nemirovski09,Shamir2013icml}. These methods are inherently
serial, in that each iteration depends on the result of the previous
iteration. Recently, parallel multicore versions of stochastic
gradient and stochastic coordinate descent have been described for
problems involving large data sets; see for example
\citep{Hogwild11nips, Richtarik12arXiv, Avron13arXiv, LiuWright13,
  Sridhar2013nips, Liu14arXivAsyRK}.

This paper proposes an asynchronous stochastic proximal
coordinate-descent algorithm, called $\AsySPCD$, for composite
objective functions.  The basic step of $\AsySPCD$, executed
repeatedly by each core of a multicore system, is as follows: Choose
an index $i \in \{1,2,\dotsc, n\}$; read $x$ from shared memory and
evaluate the $i$th element of $\nabla f$; subtract a short, constant,
positive of this partial gradient from $(x)_i$; and perform a proximal
operation on $(x)_i$ to account for the regularization term
$g_i(\cdot)$.
%
%
We use a simple model of computation that matches well to modern
multicore architectures. Each core performs its updates on centrally
stored vector $x$ in an asynchronous, uncoordinated fashion, without
any form of locking.
A consequence of this model is that the version of $x$ that is read by
a core in order to evaluate its gradient is usually not the same as the
version to which the update is made later, because $x$ is updated in the
interim by other cores.  (Generally, we denote by $\hat{x}$ the
version of $x$ that is used by a core to evaluate its component of
$\nabla f(\hat{x})$.)  We assume, however, that indefinite delays do
not occur between reading and updating: There is a bound $\tau$ such
no more than $\tau$ component-wise updates to $x$ are missed by a
core, between the time at which it reads the vector $\hat{x}$ and the
time at which it makes its update to the chosen element of $x$.  A
similar model of parallel asynchronous computation was used in
\hogwild~\citep{Hogwild11nips} and $\AsySCD$~\citep{LiuWright13}.
However, there is a key difference in this paper: {\em We do {\em not}
  assume that the evaluation vector $\hat{x}$ is a version of $x$ that
  actually existed in the shared memory at some point in time.}
Rather, we account for the fact that the components of $x$ may be
updated by multiple cores while in the process of being read by
another core, so that $\hat{x}$ may be a ``hybrid'' version that never
actually existed in memory. Our new model, which we call an
``inconsistent read'' model, is significantly closer to the reality of
asynchronous computation, and dispenses with the somewhat unsatisfying
``consistent read'' assumption of previous work. It also requires a
quite distinct style of analysis; our proofs differ substantially from
those in previous related works.

We show that, for suitable choices of steplength, our algorithm
converges at a linear rate if an ``optimal strong convexity'' property
\eqnok{eq:esc} holds. It attains sublinear convergence at a ``$1/k$''
rate for general convex functions. Our analysis also defines a
sufficient condition for near-linear speedup in the number of cores
used. This condition relates the value of delay parameter $\tau$
(which corresponds closely to the number of cores / threads used in
the computation) to the problem dimension $n$. A parameter that
quantifies the cross-coordinate interactions in $\nabla f$ also
appears in this relationship. When the Hessian of $f$ is nearly
diagonal, the minimization problem \eqref{eqn_mainproblem} is almost
separable, so higher degrees of parallelism are possible.

We review related work in
Section~\ref{sec_relatedwork}. Section~\ref{sec_alg} specifies the
proposed algorithm. Convergence results are described in
Section~\ref{sec_mainresult}, with proofs given in the
appendix. Computational experience is reported in
Section~\ref{sec_exp}. A summary and conclusions appear in
Section~\ref{sec_conclusion}.

\subsection*{Notation and Assumption}\label{sec_NA}

We use the following notation in the remainder of the paper.
\begin{itemize}
\item $\Omega$ denotes the intersection of ${\rm dom}(f)$ and ${\rm dom}(g)$
\item $S$ denotes the set on which $F$ attains its optimal value,
  which is denoted by $F^*$.
\item $\PS(\cdot)$ denotes Euclidean-norm projection onto $S$.
\item $e_i$ denotes the $i$th natural basis vector in $\mathbb{R}^n$.
\item Given a matrix $A$, we use $A_{\cdot j}$ to denote its $j$th
  column and $A_{i \cdot}$ to denote its $i$th row.
\item $\| \cdot \|$ denotes the Euclidean norm $\|\cdot \|_2$.
\item $x_j \in \mathbb{R}^n$ denotes the $j$th iterate generated by
  the algorithm.
\item $f^*_j:=f(\PS(x_j))$ and $g^*_j:=g(\PS(x_j))$.
\item $F^*:=F(\PS(x))$ denotes the optimal objective value. (Note that
  $F^*=f^*_j + g^*_j$ for any $j$.)
\item We use $(x)_i$ for the $i$th element of $x$, and $\nabla_i f(x)$
  for the $i$th element of $\nabla f(x)$.
\item Given a scalar function $h:~\R\rightarrow \R$, define the componentwise
  proximal operator
\[
\P_{i,h}(y): = {\rm arg}\min_x {1\over 2}\|x-y\|^2 +
  h((x)_i).
\]
Similarly, for the vector function $g$, we denote
\[
\P_g(y):={\rm arg}\min_x {1\over 2}\|x-y\|^2 + g(x).
\]
Note that the proximal operator is nonexpansive, that is,
$\|\P_g(x) - \P_g(y)\| \leq \|x-y\|$.
\end{itemize}

We define the following {\em optimal strong convexity} condition for a
convex function $f$ with respect to the optimal set $S$, with
parameter $l>0$:
\begin{equation}
  F(x) - F(\PS(x)) \geq {l\over 2}\|x-\PS(x)\|^2\quad \forall x\in \Omega.
  \label{eq:esc}
\end{equation}
This condition is significantly weaker than the usual strong
convexity condition; a strongly convex function $F(.)$ is an
optimally strongly convex function, but the converse is not true in
general. We provide several examples of optimally strongly convex
functions that are not strongly convex:
\begin{itemize}
\item $F(x) = {\rm constant}$.
\item $F(x) = f(Ax)$, where $f$ is a strongly convex function and
  $A$ is any matrix, possibly one with a nontrivial kernel.
\item $F(x) = f(Ax) + {\bf 1}_{X}(x)$ with strongly convex $f$, and
  arbitrary $A$, where ${\bf 1}_X(x)$ is an indicator function
  defined on a polyhedron set $X$. Note first that $y^*:=Ax^*$ is
  unique for any $x^*\in S$, from the strong convexity of $f$. The
  optimal solution set $S$ is defined by
  \[
  Ax = y^*,\quad x\in X.
  \]
  The inequality \eqref{eq:esc} clearly holds for $x \notin X$, since
  the left-hand side is infinite in this case. For $x \in X$, we have by
  the famous theorem of Hoffman \cite{Hoffman52} that there exists $c>0$ such that
  \begin{align*}
    \|Ax-y^*\|^2 = \|A(x-\PS(x))\|^2 \geq c\|x-\PS(x)\|^2.
  \end{align*}
  Then from the strong convexity of $f(x)$, we have that there exists a
  positive number $l$ such that for any $x\in X$
  \begin{align*}
    F(Ax)-F(A\PS(x)) & =  f(Ax) - f(A\PS(x)) 
    \\ & \geq
       {l \over 2}\|A(x-\PS(x))\|^2 \geq {lc \over 2}\|x-\PS(x)\|^2.
  \end{align*}
\item Squared hinge loss $F(x) = \sum_{i} \max(0, a^T_i x -
  y_i)^2$. To verify optimal strong convexity, we reformulate this 
  problem as
  \[
  \min_{t,x} \, \|t\|^2 \;\; \mbox{subject to} \; t_i \ge a_i^Tx-y_i \;\; \forall_i,
  \]
  and apply the result just derived.
\end{itemize}

Note that optimal strong convexity \eqref{eq:esc} is a weaker
version of the ``essential strong convexity'' condition used in
\citet{LiuWright13}. A concept called ``restricted strong convexity''
proposed in \citet{LaiYin13} (See Lemma~4.6) is similar in that it requires
a certain quantity to increase quadratically with distance from the
solution set, but different in that the objective is assumed to be
differentiable. Anitescu \citet{Ani00a} defines a ``quadratic growth
condition'' for (smooth) nonlinear programming in which the objective is
assumed to grow at least quadratically with distance to a local
solution in some feasible neighborhood of that solution. Since our
setting (unconstrained, nonsmooth, convex) is quite different, we
believe the use of a different term is warranted here.

Throughout this paper, we make the following assumption.
\begin{assumption} \label{ass_1}
The  solution set $S$ of \eqref{eqn_mainproblem} is nonempty.
\end{assumption}

\subsection*{Lipschitz Constants}

We define two different Lipschitz constants $\Lres$ and $\Lmax$ that
are critical to the analysis, as follows.  $\Lres$ is the {\em
  restricted Lipschitz constant} for $\nabla f$ along the coordinate
directions: For any $x\in \Omega$, for any $i=1,2,\dotsc,n$, and any
$t \in \mathbb{R}$ such that $x+te_i \in \Omega$, we have
\[
\| \nabla f(x) - \nabla f(x+te_i)\| \leq \Lres |t|.
\]
The {\em coordinate Lipschitz constant} $\Lmax$ is defined for $x$,
$i$, $t$ satisfying the same conditions as above:
\[
\|\nabla f(x) - \nabla f(x+te_i)\|_{\infty} \leq \Lmax|t|.
\]
Note that
\begin{equation} \label{eq:col}
f(x+te_i) - f(x) \leq \langle \nabla_if(x),~t \rangle + {\Lmax\over 2}t^2.
\end{equation}
We denote the ratio between these two quantities by
$\Lambda$:
\begin{equation} \label{def:Lambda}
\Lambda:={\Lres}/{\Lmax}.
\end{equation}
Making the implicit assumption that $\Lres$ and $\Lmax$ are chosen to
be the smallest values that satisfy their respective definitions, we
have from standard relationships between the $\ell_2$ and
$\ell_{\infty}$ norms that
\[
1 \le \Lambda \le \sqrt{n}.
\]

Besides bounding the nonlinearity of $f$ along various directions, the
quantities $\Lres$ and $\Lmax$ capture the interactions between the
various components in the gradient $\nabla f$. In the case of twice
continuously differentiable $f$, we can understand these interactions
by observing the diagonal and off-diagonal terms of the Hessian
$\nabla^2 f(x)$.
Let us consider upper bounds on the ratio $\Lambda$ in various
situations.  For simplicity, we suppose that $f$ is quadratic with
positive semidefinite Hessian $Q$. 
\begin{itemize}
\item
If $Q$ is sparse with at most $p$ nonzeros per row/column, we have
that
\[
\Lres = \max_i \, \| Q_{\cdot i} \|_2 \le \sqrt{p} \max_i \, \|Q_{\cdot i} \|_{\infty} = \sqrt{p} \Lmax,
\]
so that $\Lambda \le \sqrt{p}$ in this situation. 
\item If $Q$ is diagonally dominant, we have for any column $i$ that
\[
\| Q_{\cdot i} \|_2 \le Q_{ii} + \| [Q_{ji}]_{j \neq i} \|_2 \le
Q_{ii} + \sum_{j \neq i} |Q_{ji}| \le 2 Q_{ii},
\]
which, by taking the maximum of both sides, implies that $\Lambda
\le 2$ in this case.
\item Suppose that $Q=A^TA$, where $A\in\mathbb{R}^{m\times n}$ is a
  random matrix whose entries are i.i.d from $\mathcal{N}(0,1)$. (For
  example, $f$ could be the linear least-squares objective
  $f(x)=\frac12 \|Ax-b\|^2$.)  We show in \cite{LiuWright13} that
  $\Lambda$ is upper-bounded roughly by $1+\sqrt{n/m}$ in this case.


\end{itemize}

\section{Related Work}\label{sec_relatedwork}

We have surveyed related work on coordinate descent and stochastic
gradient methods in a recent report \cite{LiuWright13}. Our discussion
there included non-stochastic, cyclic coordinate descent methods
\citet{Tseng01,LuoTseng92,WangLin13,Beck13,TseY06,TseY07a, Saha13},
synchronous parallel methods that distribute the work of function and
gradient evaluation
\cite{Ferris94,Mangasarian95,Ma12,Boyd11,Duchi12,AgarwalD12,Cotter11,Shalev-Shwartz2013},
and asynchronous parallel stochastic gradient methods (including the
randomized Kaczmarz algorithm)
\cite{Hogwild11nips,Liu14arXivAsyRK}. We make some additional comments
here on related topics, and include some recent references from this
active research area.


{\em Stochastic coordinate descent} can be viewed as a special case of
stochastic gradient, so analysis of the latter approach can be
applied, to obtain for example a sublinear $1/k$ rate of convergence
in expectation for strongly convex functions; see, for example
\cite{Nemirovski09}.  However, stochastic coordinate descent is
``special'' in that it is possible to guarantee improvement in the
objective at every step. Nesterov \citet{Nesterov12} studied the
convergence rate for a stochastic block coordinate descent method for
unconstrained and separably constrained convex smooth optimization,
proving linear convergence for the strongly convex case and a
sublinear $1/k$ rate for the convex case. Richt{\'a}rik and
Tak{\'a}{\v c} \citet{Richtarik11} and Lu and Xiao \citet{LuXiao13}
extended this work to composite minimization, in which the objective
is the sum of a smooth convex function and a separable nonsmooth
convex function, and obtained similar (slightly stronger) convergence
results. Stochastic coordinate descent is extended by Necoara and
Patrascu \citet{Necoara13a} to convex optimization with a single
linear constraint, randomly updating \emph{two} coordinates at a time
to maintain feasibility.

In the class of {\em synchronous parallel methods} for coordinate
descent, Richt{\'a}rik and Tak{\'a}{\v c} \citet{Richtarik12arXiv}
studied a synchronized parallel block (or minibatch) coordinate
descent algorithm for composite optimization problems of the form
\eqref{eqn_mainproblem}, with a block separable regularizer $g$. At
each iteration, processors update the randomly selected coordinates
concurrently and synchronously. Speedup depends on the sparsity of the
data matrix that defines the loss functions. A similar synchronous
parallel method was studied in \citet{Necoara13b} and
\citet{Bradley11}; the latter focuses on the case of $g(x)=\|x\|_1$.
Scherrer et al.  \citet{ScherrerTHH12} make greedy choices of multiple
blocks of variables to update in parallel. Another greedy way of
selecting coordinates was considered by Peng et al. \citet{Yin13}, who
also describe a parallel implementation of FISTA, an accelerated
first-order algorithm due to Beck and Teboulle~\citet{BeckT09}. Fercoq
and Richt{\'a}rik~\citet{Fercoq13a} consider a variant of
\eqnok{eqn_mainproblem} in which $f$ is allowed to be nonsmooth. They
apply Nesterov's smoothing scheme to obtain a smoothed version and
update multiple blocks of coordinates using block coordinate descent
in parallel. Sublinear convergence rate is established for both
strongly convex and weakly convex cases. 
%
%
Fercoq and Richt{\'a}rik \citet{Fercoq13b} proposed a variant of
Nesterov's accelerated scheme to accelerate the synchronous parallel
block coordinate algorithm of \citet{Richtarik12arXiv}, proving an
improved sublinear convergence rate for weakly convex problems. This
variant avoids the disadvantage of the original Nesterov acceleration
scheme \citet{Nesterov12}, which requires $O(n)$ complexity per
iteration, even on sparse data.
Facchinei, Sagratella, and Scutari \cite{FacSS13}
  consider a general framework for synchronous block coordinate
  descent methods with separable regularizers, in which the block
  subproblems may be solved inexactly. However, the block to be
  updated at each step is not chosen randomly; it must contain a
  component that is furthest from optimality, in some sense.

We turn now to {\em asynchronous parallel methods}. Bertsekas and
Tsitsiklis \citet{Bertsekas89} described an asynchronous method for
fixed-point problems $x = q(x)$ over a separable convex closed
feasible region. (The optimization problem \eqnok{eqn_mainproblem} can
be formulated in this way by defining $q(x) := \mathcal{P}_{\alpha
  g}[(I - \alpha \nabla f)(x)]$ for a fixed $\alpha>0$.) They use an
inconsistent-read model of asynchronous computation, and establish
linear convergence provided that components are not neglected
indefinitely and that the iteration $x=q(x)$ is a maximum-norm
contraction. The latter condition is quite strong. In the case of $g$
null and $f$ convex quadratic in \eqnok{eqn_mainproblem} for instance,
it requires the Hessian to satisfy a diagonal dominance condition ---
a stronger condition than strong convexity. By comparison, $\AsySCD$
\citep{LiuWright13} guarantees linear convergence under an ``essential
strong convexity'' condition, though it assumes a consistent-read
model of asynchronous computation. Elsner et al.~\citet{Elsner92}
considered the same fixed point problem and architecture as
\citet{Bertsekas89}, and describe a similar scheme. Their scheme
appears to require locking of the shared-memory data structure for $x$
to ensure consistent reading and writing.  Frommer and
Szyld~\citet{Frommer00} give a comprehensive survey of asynchronous
methods for solving fixed-point problems.

Liu et al. \citet{LiuWright13} followed the asynchronous
consistent-read model of $\hogwild$~to develop an asynchronous
stochastic coordinate descent ($\AsySCD$) algorithm and proved
sublinear ($1/k$) convergence on general convex functions and a linear
convergence rate on functions that satisfy an ``essential strong
convexity'' property. Sridhar et al.~\citet{Sridhar2013nips} developed
an efficient LP solver by relaxing an LP problem into a
bound-constrained QP problem, which is then solved by $\AsySCD$.

Liu et al.~\citet{Liu14arXivAsyRK} developed an asynchronous parallel
variant of the randomized Kaczmarz algorithm for solving a general
consistent linear system $Ax=b$, proving a linear convergence rate.
Avron et al. \citet{Avron13arXiv} proposed an asynchronous solver for
the system $Qx=c$ where $Q$ is a symmetric positive definite matrix,
proving a linear convergence rate.  This method is essentially an
asynchronous stochastic coordinate descent method applied to the
strongly convex quadratic optimization problem $\min_x \, {1\over
  2}x^TQx-c^Tx$. The paper considers both inconsistent- and
consistent-read cases are considered, with slightly different
convergence results.

\section{Algorithm}\label{sec_alg}

\begin{figure}[t!]
\center
\vspace{-5mm}
\includegraphics[scale=0.5]{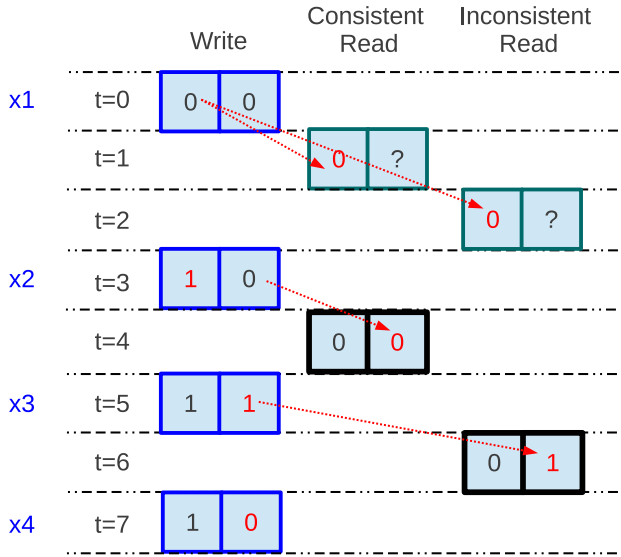}
\vspace{-20mm}
\caption{Time sequence of writes and reads of a two-variable vector,
  showing instances of consistent and inconsistent reading. The left
  column shows the initial vector at time $0$, stored in shared
  memory, with updates to single components at times 3, 5, and 7. The
  middle column shows a consistent read, in which the first component
  is read at time 1 and the second component is read at time 4.  The
  read vector is equal to the shared-memory vector at time 0.  The
  right column shows an inconsistent read, in which the first
  component is read at time 2 and the second component is read at time
  6. Because of intervening writes to these components, the read vector
  does not match the versions that appeared in shared memory at any
  time point.}
\label{fig_inconsistentread}
\end{figure}

In our algorithm \AsySPCD, multiple processors have access to a shared
data structure for the vector $x$, and each processor is able to
compute a randomly chosen element of the gradient vector $\nabla
f(x)$. Each processor repeatedly runs the following proximal
coordinate descent process. (Choice of the steplength parameter
$\gamma$ is discussed further in the next section.)
\begin{itemize}
\item[R:] Choose an index $i \in \{1,2,\dotsc,n\}$ at random, read $x$
  into the local storage location $\hx$, and evaluate $\nabla_i
  f(\hat{x})$;
\item[U:] Update component $i$ of the shared $x$ by taking a step of
  length $\gamma/\Lmax$ in the direction $-\nabla_i f(\hx)$, follows
  by a proximal operation defined as follows:\footnote{Our
    analysis assumes that no other process modifies $x_i$ while this
    proximal operation is being computed. As we explain in Section~\ref{sec_exp},
    our practical implementation actually assigns each coordinate
    $x_i$ to a single core, and allows only that core to update $x_i$,
    so this issue does not arise. An alternative implementation,
    pointed out by a referee, would be to use a ``compare-and-swap''
    atomic instruction to implement the update. This operation would
    perform the update only if $x_i$ was not changed while the update
    was being computed.}
\[
x \leftarrow \P_{i,{\gamma\over
    \Lmax}g_{i}}\left(x-\frac{\gamma}{\Lmax}
e_{i}\nabla_{i}f(\hx)\right).
\]
\end{itemize}
Notice that each step changes just a single element of $x$, that is,
the $i$th element. Unlike standard proximal coordinate descent, the
value $\hx$ at which the coordinate gradient is calculated usually
differs from the value of $x$ to which the update is applied, because
while the processor is evaluating its gradient, other processors may
repeatedly update the value of $x$ stored in memory.  As mentioned
above, we use an ``inconsistent read'' model of asynchronous
computation here, in contrast to the ``consistent read'' models of
$\AsySCD$ \citep{LiuWright13} and $\hogwild$ \citep{Hogwild11nips}.
Figure~\ref{fig_inconsistentread} shows how inconsistent reading can
occur, as a result of updating of components of $x$ while it is being
read.  Consistent reading can be guaranteed by means of a
software lock, but such a mechanism degrades parallel performance
significantly. In fact, the implementations of $\hogwild$ and
$\AsySCD$ described in the papers \citep{Hogwild11nips,LiuWright13} do
not use any software lock, and in this respect the computations in
those papers are not quite compatible with their analysis.

\begin{algorithm}
\caption{Asynchronous Stochastic Coordinate Descent Algorithm $x_{J}=\AsySPCD(x_0, \gamma, J)$}
\label{alg_ascd}
\begin{algorithmic}[1]
\REQUIRE $x_0$, $\gamma$, and $J$
\ENSURE $x_{J}$
\STATE Initialize $j \leftarrow 0$;
\WHILE{$j < J$}
\STATE Choose $i(j)$ from $\{1,2,\dotsc,n\}$ with equal probability;
\STATE $x_{j+1} \leftarrow \P_{i(j),{\gamma\over \Lmax}g_{i(j)}}\left(x_j-\frac{\gamma}{\Lmax} e_{i(j)} \nabla_{i(j)}f(\hx_{j})\right)$; \label{step_proj}
\STATE $j \leftarrow  j +1$;
\ENDWHILE
\end{algorithmic}
\end{algorithm}

The ``global'' view of algorithm $\AsySPCD$ is shown in {\bf
  Algorithm} \ref{alg_ascd}.  To obtain this version from the
``local'' version, we introduce a counter $j$ to track the total
number of updates applied to $x$, so that $x_j$ is the state of $x$ in
memory after update $j$ is performed. We use $i(j)$ to denote the
component that is updated at iteration $j$, and $\hx_j$ for value of
$x$ that is used in the calculation of the gradient element $\nabla
f_{i(j)}$.  The components of $\hx_{j}$ may have different ages. Some
components may be current at iteration $j$, others may not reflect
recent updates made by other processors. We assume however that there
is an upper bound of $\tau$ on the age of each component, measured in
terms of updates. $K(j)$ defines an iterate set such that
\[
x_j = \hx_{j} + \sum_{d\in K(j)} (x_{d+1}-x_d).
\]
One can see that $d\leq j-1$, $\forall d\in K(j)$. Here we assume $\tau$
to be the upper bound on the age of all elements in $K(j)$, for all
$j$, so that $\tau \geq j-\min\{d~|~d\in K(j)\}$. We assume further
that $K(j)$ is ordered from oldest to newest index (that is, smallest
to largest). Note that $K(j)$ is empty if $x_j=\hx_j$, that is, if the
step is simply an ordinary stochastic coordinate gradient update. The
value of $\tau$ corresponds closely to the number of cores involved in
the computation provided that computation of the update for each
  component of $x$ costs roughly the same.

\section{Main Results}\label{sec_mainresult}

This section presents results on convergence of $\AsySPCD$. The
theorem encompasses both the linear rate for optimally strongly convex
$f$ and the sublinear rate for general convex $f$. The result depends
strongly on the delay parameter $\tau$. The proofs are highly
technical, and are relegated to Appendix~\ref{app:con}. We note the
proof techniques differ significantly from those used for the
consistent-read algorithms of \citet{Hogwild11nips} and
\citet{LiuWright13}.

We start by describing the key idea of the algorithm, which is
reflected in the way that it chooses the steplength parameter
$\gamma$.  Denoting $\bar{x}_{j+1}$ by
\begin{equation}
\bar{x}_{j+1} :=\mathcal{P}_{{\gamma \over \Lmax}g}\left(x_j-{\gamma\over \Lmax} \nabla f(\hx_{j})\right),
\label{eqn_xj1}
\end{equation}
we can see that 
\begin{equation} \label{xbarij}
(x_{j+1})_{i(j)} = (\bar{x}_{j+1})_{i(j)}, \quad (x_{j+1})_i = (x_j)_i
\;\; \mbox{for $i \neq i(j)$},
\end{equation}
so that $x_{j+1}-x_j = [(\bar{x}_{j+1})_{i(j)} -(x_j)_{i(j)} ] e_{i(j)}$.  Thus, we
have
\[ 
\E_{i(j)}(x_{j+1} - x_j) = 
\frac{1}{n} \sum_{i=1}^n [ (\bar{x}_{j+1})_{i} - (x_j)_i ] e_i =
\frac{1}{n} [ \bar{x}_{j+1}-x_j].
\] 
Therefore, we can view $\bar{x}_{j+1}-x_j$ as capturing the expected
behavior of $x_{j+1}-x_j$. Note that when $g(x)=0$, we have
$\bar{x}_{j+1}-x_j=-({\gamma}/{\Lmax})\nabla f(\hx_j)$, a standard
negative-gradient step. The choice of steplength parameter $\gamma$
entails a tradeoff: We would like $\gamma$ to be long enough that
significant progress is made at each step, but not so long that the
gradient information computed at $\hx_j$ is stale and irrelevant by
the time the update is applied to $x_j$. We enforce this tradeoff by
means of a bound on the ratio of expected squared norms on
$x_j-\bar{x}_{j+1}$ at successive iterates; specifically,
\begin{equation} \label{eq:ratio_unc}
\E\|x_{j-1}-\bar{x}_j\|^2 \leq \rho \E\|x_j - \bar{x}_{j+1}\|^2,
\end{equation}
where $\rho > 1$ is a user defined parameter.  The analysis becomes a
delicate balancing act in the choice of $\rho$ and steplength $\gamma$
between aggression and excessive conservatism. We find, however, that
these values can be chosen to ensure steady convergence for the
asynchronous method at a linear rate, with rate constants that are
almost consistent with a standard short-step proximal full-gradient
descent, when the optimal strong convexity condition \eqnok{eq:esc} is
satisfied.

Our main convergence result is the following.

\begin{theorem} \label{thm_2}
Suppose that Assumption~\ref{ass_1} is satisfied. Let $\rho$ be a
constant that satisfies $\rho>1+4/\sqrt{n}$, and define the quantities
$\theta$, $\theta'$, and $\psi$ as follows:
\beq \label{eq:defpsic}
\theta:=\frac{\rho^{(\tau+1)/2}-\rho^{1/2}}{\rho^{1/2}-1}, \quad \theta':=\frac{\rho^{(\tau+1)}-\rho}{\rho-1}, \quad \psi:= 1 + \frac{\tau\theta'}{n} + \frac{2\Lambda\theta}{\sqrt{n}}.  
\eeq 
Suppose that the steplength parameter $\gamma>0$ satisfies the following two bounds:
  \beq \label{eq:boundgammac}
\gamma\le {1\over \psi}, \quad\gamma \le \frac{\sqrt{n}(1-\rho^{-1})-4}{4(1+\theta)\Lambda}.
\eeq
Then we have
\begin{equation}
\E\|x_{j-1}-\bar{x}_j\|^2 \leq \rho \E \|x_j-\bar{x}_{j+1}\|^2,
\quad j=1,2,\dotsc. \label{eqn_thm2_1}
\end{equation}
If the optimal strong convexity property~\eqref{eq:esc}
holds with $l>0$, we have for $j=1,2,\dotsc$ that
\begin{align}
\nonumber
&\E \|x_{j} - \PS(x_{j})\|^2 + \frac{2\gamma}{\Lmax}(\E F(x_{j}) - F^*) 
\\ &\quad 
\leq \left(1- \frac{l} {n(l+\gamma^{-1}\Lmax)} \right)^{j} \left(\|x_0-\PS(x_0)\|^2 + \frac{2\gamma}{\Lmax}(F(x_0) - F^*)\right), 
\label{eqn_thm2_3}
\end{align}
while for general smooth convex function $f$, we have
\beq  \label{eqn_thm2_2}
\E F(x_{j})- F^* \leq \frac{n(\|x_0-\PS(x_0)\|^2\Lmax+ 2\gamma(F(x_0)- F^*))}{2\gamma (n+j)}.
\eeq
\end{theorem}

The following corollary proposes an interesting particular choice for
the parameters for which the convergence expressions become more
comprehensible. The result requires a condition on the delay bound
$\tau$ in terms of $n$ and the ratio $\Lambda$.

\begin{corollary} \label{co:thm_2}
Suppose that Assumption~\ref{ass_1} holds and that
\beq \label{eq:boundtauc}
4e\Lambda(\tau + 1)^2 \leq \sqrt{n}.
\eeq
If we choose
\beq \label{eq:choicerhoc}
\rho=\left(1+{4e\Lambda(\tau+{1}) \over  \sqrt{n}}\right)^2,
\eeq
then the steplength $\gamma=1/2$ will satisfy the bounds
\eqnok{eq:boundgammac}. In addition, when the optimal strong convexity
property \eqnok{eq:esc} holds with $l>0$, we have for $j=1,2,\dotsc$
that
\begin{equation}
 \label{eqn_thm_2_good_c}
\E F(x_j)-F^*
 \le \left( 1-\frac{l}{n(l+2\Lmax)} \right)^j (\Lmax \|x_0-\PS(x_0)\|^2 + F(x_0)-F^*),
\end{equation}
while for the case of general convex $f$, we have
\beq \label{eqn_thm_3_good_c}
\E F(x_j)-F^* \le  \frac{n(\Lmax \|x_0-\PS(x_0)\|^2 + F(x_0)-F^*)}{j+n}.
\eeq
\end{corollary}

We note that the linear rate \eqnok{eqn_thm_2_good_c} is broadly
consistent with the linear rate for the classical steepest descent
method applied to strongly convex functions, which has a rate constant
of $(1-2l/L)$, where $L$ is the standard Lipschitz constant for
$\nabla f$. Suppose we assume (not unreasonably) that $n$ steps of
stochastic coordinate descent cost roughly the same as one step of
steepest descent, and that $l \le \Lmax$. It follows from
\eqnok{eqn_thm_2_good_c} that $n$ steps of stochastic coordinate
descent would achieve a reduction factor of about
\[
1-\frac{l}{2 \Lmax + l} \leq 1-\frac{l}{3\Lmax}, 
\]
\jlcommentsolved{(we need assume $\Lmax \ge l$ here, because sometimes
  $l$ could larger than $\Lmax$.)} so a standard argument would
suggest that stochastic coordinate descent would require about $6
\Lmax/L$ times more computation. Since $\Lmax/L \in [1/n,1]$, the
stochastic asynchronous approach may actually require less
computation.  It may also gain an advantage from the parallel
asynchronous implementation. A parallel implementation of standard
gradient descent would require synchronization and careful division of
the work of evaluating $\nabla f$, whereas the stochastic approach can
be implemented in an asynchronous fashion.

For the general convex case, \eqnok{eqn_thm_3_good_c} defines a
sublinear rate, whose relationship with the rate of standard gradient
descent for general convex optimization is similar to the previous
paragraph.

Note that the results in Theorem~\ref{thm_2} and
Corollary~\ref{co:thm_2} are consistent with the analysis for
constrained $\AsySCD$ in \citet{LiuWright13}, but this paper considers
the more general case of composite optimization and the
inconsistent-read model of parallel computation.

As noted in Section~\ref{sec:intro}, the parameter $\tau$ corresponds
closely to the number of cores that can be involved in the
computation, since if all cores are working at the same rate, we would
expect each other core to make one update between the times at which
$x$ is read and (later) updated.  If $\tau$ is small enough that
\eqref{eq:boundtauc} holds, the analysis indicates that near-linear
speedup in the number of processors is achievable. A small value for
the ratio $\Lambda$ (not much greater than $1$) implies a greater
degree of potential parallelism. As we note at the end of
Section~\ref{sec:intro}, this ratio tends to closer to $1$ than to
$\sqrt{n}$ in some important applications. In these situations, the
bound \eqref{eq:boundtauc} indicates that $\tau$ can vary like
$n^{1/4}$ without affecting the iteration-wise convergence rate, and 
yielding near-linear speedup in the number of cores.
 This quantity is
consistent with the analysis for constrained $\AsySCD$ in
\citet{LiuWright13} but weaker than the unconstrained $\AsySCD$ (which
allows the maximal number of cores being $O(n^{1/2})$). A further comparison 
is with the asynchronous randomized Kaczmarz algorithm \citep{Liu14arXivAsyRK} 
which allows $O(m)$ cores to be used efficiently when solving a consistent sparse linear system. 

We conclude this section with a high-probability bound. The result
follows immediately from Markov's inequality. See Theorem~3 in
\citet{LiuWright13} for a related result and complete proof.
\begin{theorem} \label{thm:co2}
Suppose that the conditions of Corollary~\ref{co:thm_2} hold,
including the choice of $\rho$. Then for $\epsilon>0$ and
$\eta\in(0,1)$, we have that 
\begin{equation} \label{eq:thm:co2_2}
\mathbb{P}\left(F(x_j)-F^*\le \epsilon \right) \ge 1-\eta,
\end{equation}
provided that one of the following conditions holds.  In the
optimally strongly convex case~\eqnok{eq:esc} with $l>0$, we require
\[
j \geq \frac{n(l+2\Lmax)}{l} \left|{\log{\Lmax \|x_0-\PS(x_0)\|^2 + F(x_0)-F^*\over \epsilon\eta}}\right|, 
\]
iterations, while in the general convex case, it suffices that
\[
j \geq \frac{n(\Lmax \|x_0-\PS(x_0)\|^2+F(x_0)-F^*)}{\epsilon\eta} - n.
\]
\end{theorem}

\section{Experiments} \label{sec_exp}

This section presents some results to illustrate the effectiveness of
$\AsySPCD$, in particular, the fact that near-linear speedup can be
observed on a multicore machine. We note that more comprehensive
experiments can be found in \citet{LiuWright13} and
\citet{Sridhar2013nips}, for unconstrained and box-constrained
problems.  Although the analysis in \citet{LiuWright13} assumes
consistent read, it is not enforced in the implementation, so apart
from the fact that we now include a prox-step to account for the
regularization term, the implementations in \citet{LiuWright13} and
\citet{Sridhar2013nips} are quite similar to the one employed in this
section.

We apply our code for $\AsySPCD$ to the following
``$\ell_2$-$\ell_1$'' problem:
\[
\min_x\quad {1\over 2}\|Ax-b\|^2 + \lambda\|x\|_1 \equiv {1\over
  2}x^TA^TAx - b^TAx + {1\over 2}b^Tb + \lambda \|x\|_1.
\]
The elements of $A\in \R^{m\times n}$ are selected i.i.d. from a
Gaussian $\mathcal{N}(0, 1)$ distribution. To construct a sparse true
solution $x^*\in \R^n$, given the dimension $n$ and sparsity $s$, we
select $s$ entries of $x^*$ at random to be nonzero \sjwcommentsolved{I am
  totally not understanding this. Do you mean to say that $x^*$ is
  sparse with $s$ nonzeros? $A$ is dense, right?}  and
$\mathcal{N}(0,1)$ normally distributed, and set the rest to zero. The
measurement vector $b\in \R^{m}$ is obtained by $b = Ax^* + \epsilon$,
\sjwcommentsolved{Do you mean $x^*$ here?}  where elements of the noise
vector $\epsilon\in \R^m$ are i.i.d.  $\mathcal{N}(0, \sigma^2)$,
where the value of $\sigma$ controls the signal-to-noise ratio.

Our experiments run on $1$ to $10$ threads of an Intel Xeon machine,
with all threads sharing a single memory socket. Our implementations
deviate modestly from the version of $\AsySPCD$ described in
Section~\ref{sec_alg}. We compute $Q:=A^TA\in\R^{n\times n}$ and
$c:=A^Tb\in \R^n$ offline. $Q$ and $c$ are partitioned into slices
(row submatrices) and subvectors (respectively) of equal size, and
each thread is assigned one submatrix from $Q$ and the corresponding
subvector from $c$. During the algorithm, each thread updates the
elements of $x$ corresponding to its slice of $Q$, in order. After one
scan, or ``epoch'' is complete, it reorders the indices randomly, then
repeats the process.  This scheme essentially changes the scheme from
sampling with replacement (as analyzed) to sampling without
replacement, which has demonstrated empirically better performance on
many related problems. (The same advantage is noted in the
implementations of $\hogwild$ \citep{Hogwild11nips}.)

\begin{figure}[htp]
  \centering
 \includegraphics[width=0.45\textwidth]{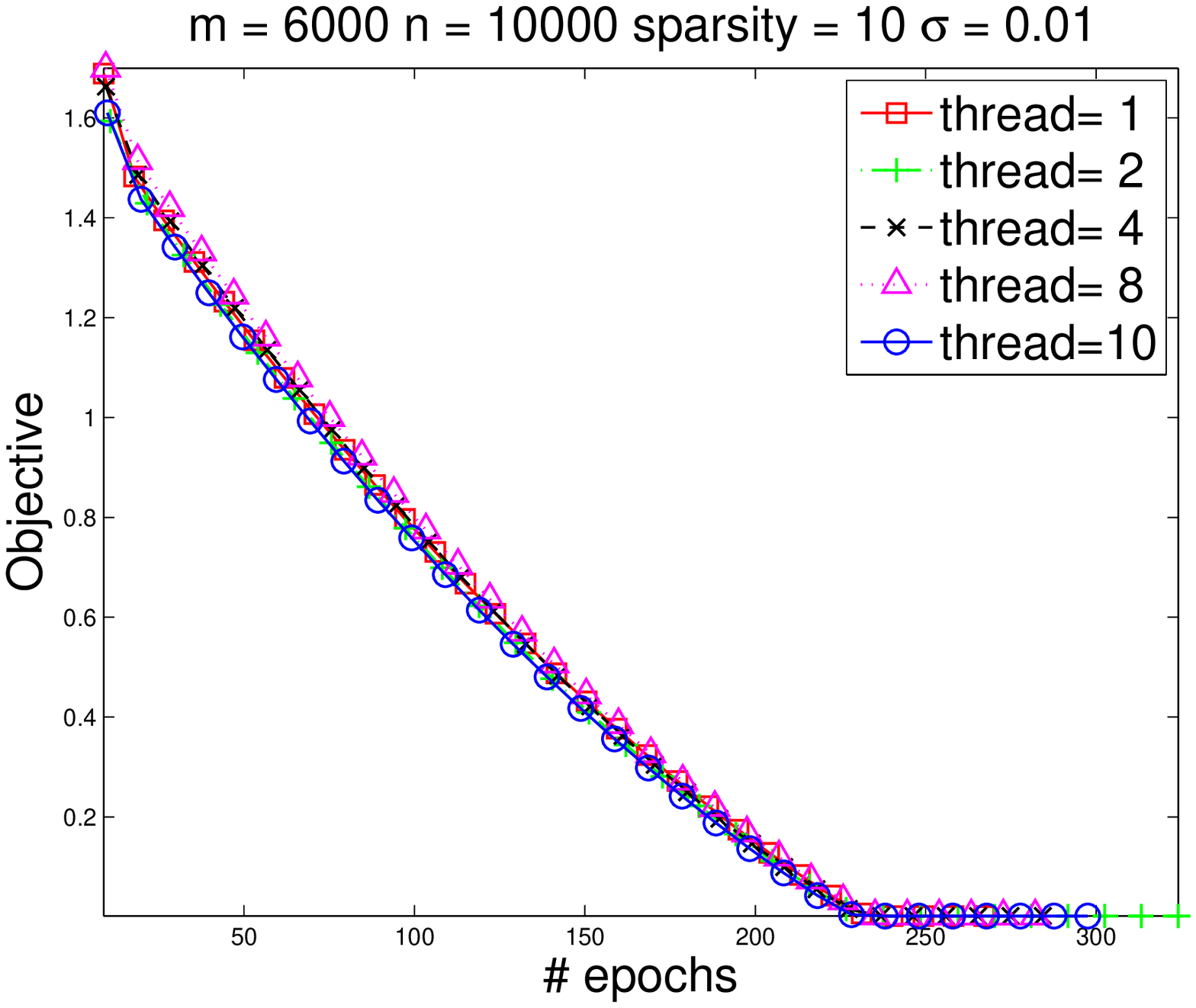} \;\;
 \includegraphics[width=0.45\textwidth]{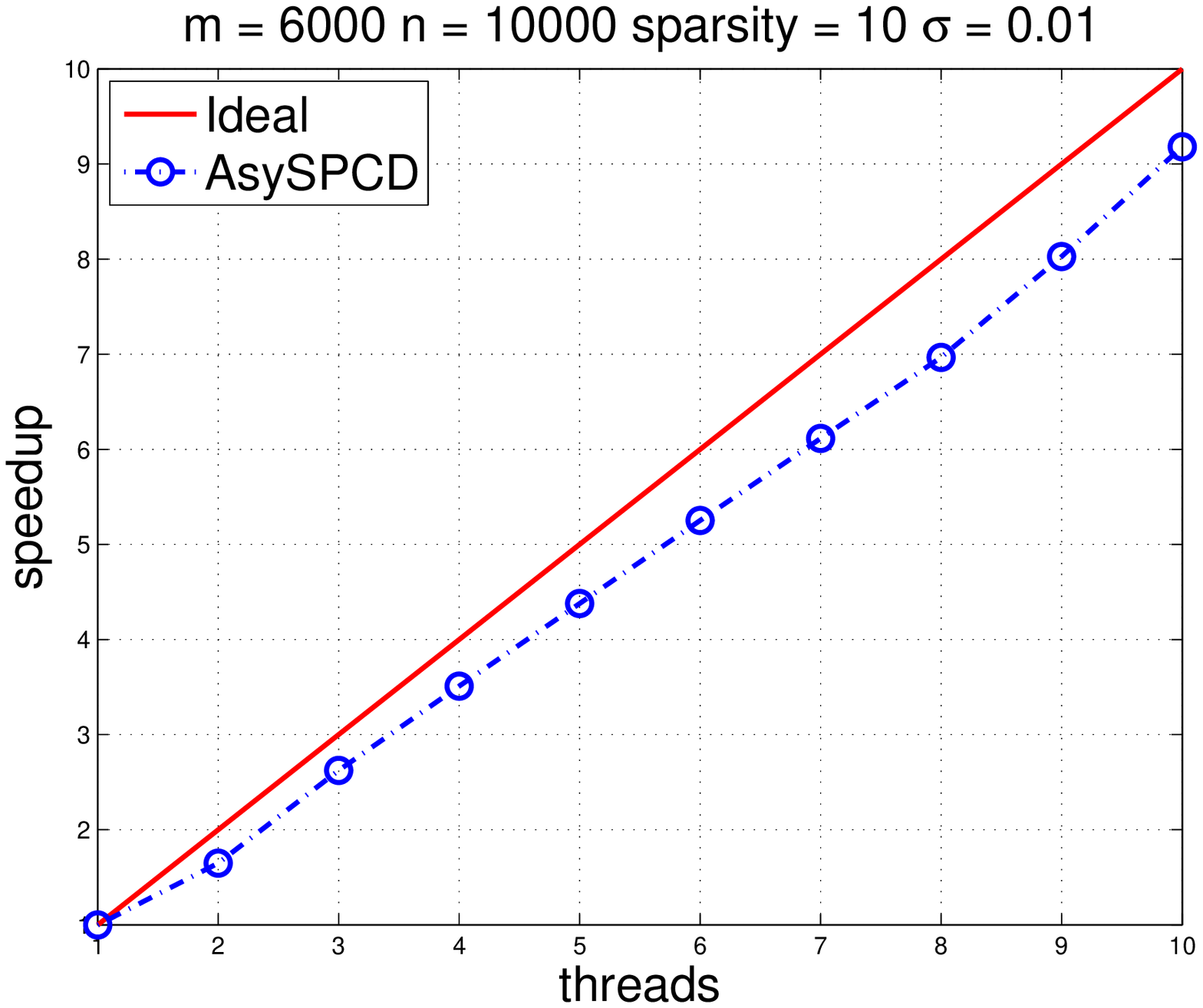}
\caption{The left graph plots objective function vs epochs for 1, 2,
  4, 8, and 10 cores. The right graph shows speedup obtained for
  implementation on 1-10 cores, plotted against the ideal (linear)
  speedup.}
\label{fig:LASSO_1}
\end{figure}

\begin{figure}[htp]
  \centering
 \includegraphics[width=0.45\textwidth]{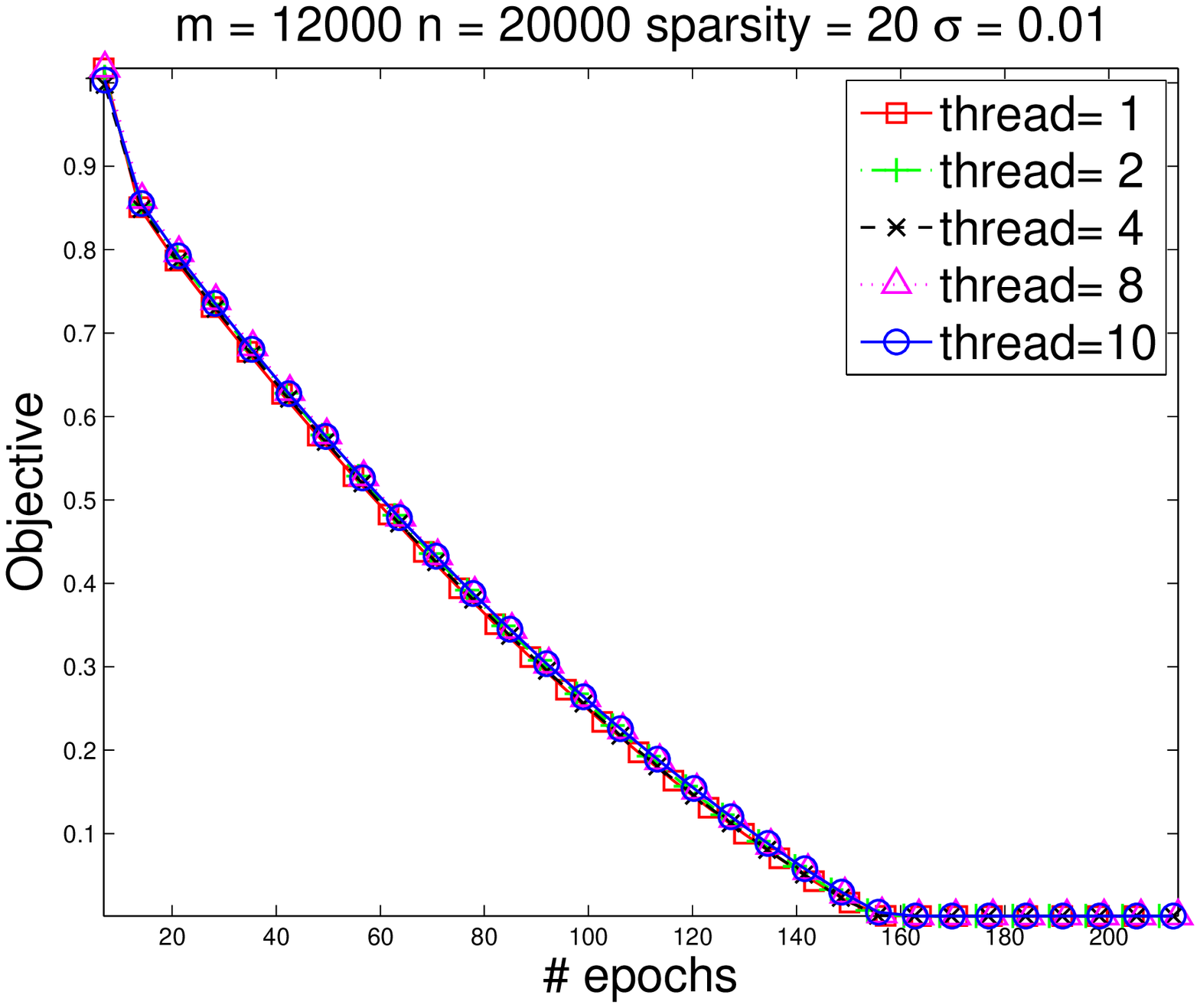} \;\;
 \includegraphics[width=0.45\textwidth]{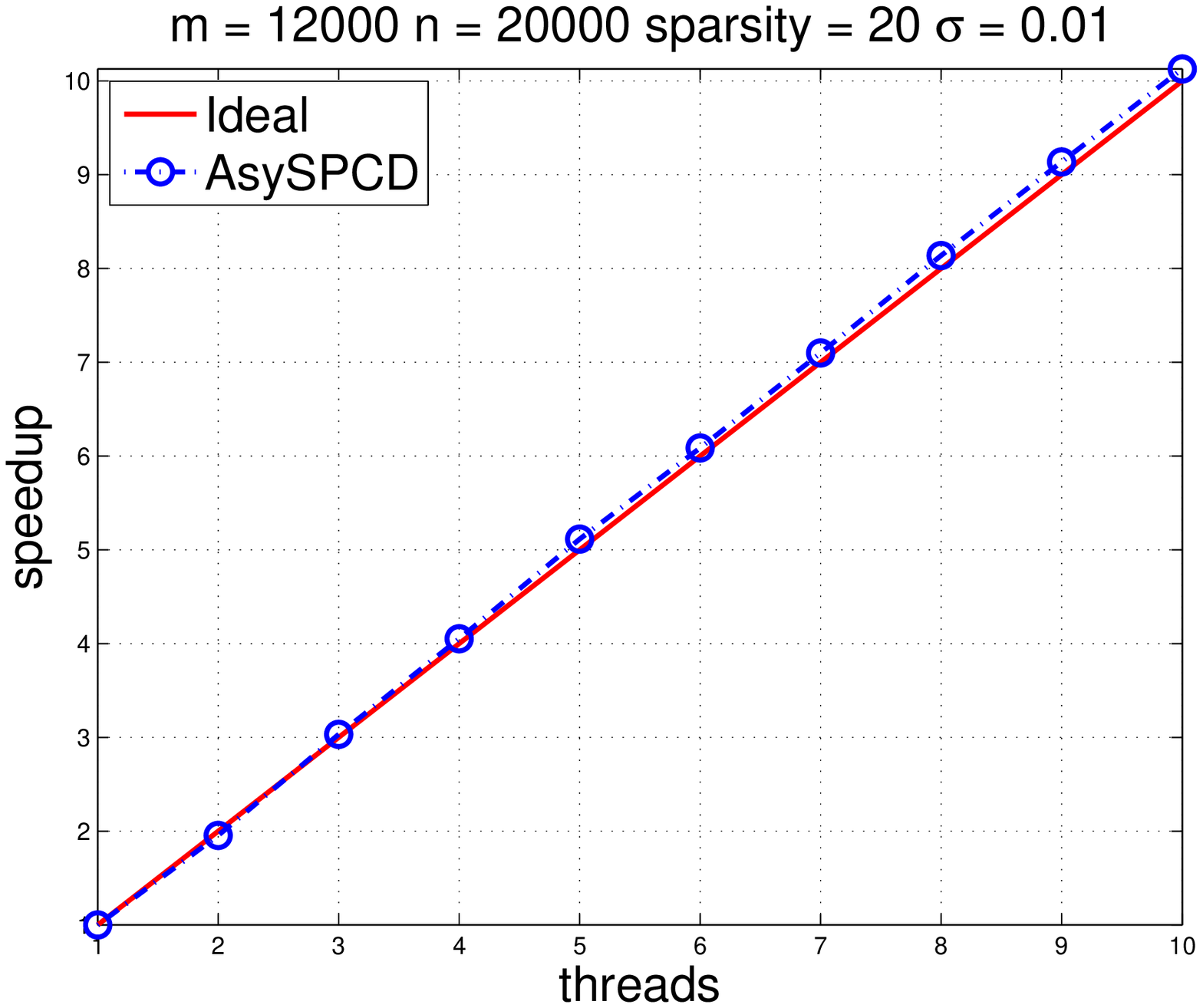}
\caption{The left graph plots objective function vs epochs for 1, 2,
  4, 8, and 10 cores. The right graph shows speedup obtained for
  implementation on 1-10 cores, plotted against the ideal (linear)
  speedup.}
\label{fig:LASSO_2}
\end{figure}

We choose $\sigma = 0.01$ with $m=6000$, $n=10000$, and $s=10$ in
Figure~\ref{fig:LASSO_1} and $m=12000$, $n=20000$, and $s=20$ in
Figure~\ref{fig:LASSO_2}. We set $\lambda = 20\sqrt{m\log (n)}\sigma$
(a value of the order of $\sqrt{m\log (n)}\sigma$ is suggested by
compressed sensing theory) and the steplength $\gamma$ is set as $1$
in both figures. In both cases, we can estimate the ratio $\Lambda =
\Lres/\Lmax$ roughly by $1+\sqrt{n/m} \approx 2.3$, as suggested at
the end of Section~\ref{sec:intro}.  Our final computed values of $x$
have nonzeros in the same locations as the chosen solution $x^*$,
though the values differ, because of the noise in $b$.

The left-hand graph in each figure indicates the number of threads /
cores and plots objective function value vs epoch count, where one
epoch is equivalent to $n$ iterations. \sjwcommentsolved{What about
  the quality of solution? If you designed $x^*$ with $s$ nonzeros,
  does the computer solution have close to $s$ nonzeros?} Note that
the curves are almost overlaid, indicating that the total workload
required for $\AsySPCD$ is almost independent of the number of cores
used in the computation. This observation validates our result in
Corollary~\ref{co:thm_2}, which indicates that provided $\tau$ is
below a certain threshold, it does not seriously affect the rate of
convergence, as a function of total computation performed.  The
right-hand graph in each figure shows speedup when executed on
different numbers of cores. Near-linear speedup is observed in
Figure~\ref{fig:LASSO_2}, while there is a slight dropoff for the
larger numbers of cores in Figure~\ref{fig:LASSO_1}. The difference
can be explain by the smaller dimension of the problem illustrated in
Figure~\ref{fig:LASSO_1}. Referring to our threshold value
  \eqref{eq:boundtauc} that indicates dimensions above which linear
  speedup should be expected, we have by setting $\Lambda \approx 2.3$
  (as discussed above) and $\tau=10$ (the maximum number of threads
  used in this experiment) that the left-hand side of
  \eqref{eq:boundtauc} is approximately 3000, while the right-hand
  side is $100$ (for Figure~\ref{fig:LASSO_1}) and approximately $141$
  (for Figure~\ref{fig:LASSO_1}). As expected, our analysis is quite
  conservative; near-linear speedup is observed even when the
  threshold \eqref{eq:boundtauc} is violated significantly.

\section{Conclusions} \label{sec_conclusion}

This paper proposes an asynchronous parallel proximal stochastic
coordinate descent algorithm for minimizing composite objectives of
the form \eqnok{eqn_mainproblem}. Sublinear convergence (at rate
$1/k$) is proved for general convex functions, with stronger linear
convergence results for problems that satisfy the optimal strong
convexity property \eqnok{eq:esc}. Our analysis indicates the extent
to which parallel implementations can be expected to yield near-linear
speedup, in terms of a parameter that quantifies the cross-coordinate
interactions in the gradient $\nabla f$ and a parameter $\tau$ that
bounds the delay in updating. Our computational experience confirms
that the linear speedup properties suggested by the analysis can be
observed in practice.

\section*{Acknowledgments}

The authors thank the editor and both referees for their valuable
comments. Special thanks to Dr. Yijun Huang for her implementation of
$\AsySPCD$, which was used here to obtain computational results.


\appendix

\section{Proofs of Main Results} \label{app:con}

This section provides the proofs for the main convergence results.  We
start with some preliminaries, then proceed to proofs of
Theorem~\ref{thm_2} and Corollary~\ref{co:thm_2}.

\subsection{Preliminaries}

Note that the component indices $i(0), i(1),\dotsc,i(j), \dotsc$ in
Algorithm~\ref{alg_ascd} are independent random variables. We use $\E$
to denote the expectation over all random variables, and $\E_{i(j)}$
to denote the conditional expectation in term of $i(j)$ given $i(0),
i(1), \dotsc, i(j-1)$. We also denote \beq\label{eq:defdelta}
(\Delta_j)_{i(j)}:=(x_j-x_{j+1})_{i(j)},
\eeq
and formulate the update in Step~\ref{step_proj} of
Algorithm~\ref{alg_ascd} in the following way:
\[
x_{j+1}=\arg\min_{x} \, \langle \nabla_{i(j)}f(\hx_{j}), (x-x_j)_{i(j)} \rangle + \frac{\Lmax}{2\gamma}\|x-x_j\|^2+ g_{i(j)}((x)_{i(j)}).
\]
(Note that $(x_{j+1})_i=(x_j)_i$ for $i\neq i(j)$.) From the optimality condition for this formulation (see (41) in \citep{Tseng08}), we have for all $x$ that
\begin{align*}
\left\langle (x-x_{j+1})_{i(j)}, \nabla_{i(j)}f(\hx_{j}) - {\Lmax\over \gamma}(\Delta_j)_{i(j)} \right\rangle + g_{i(j)}((x)_{i(j)}) - g_{i(j)}((x_{j+1})_{i(j)}) \geq 0.
\end{align*}
By rearranging this expression and substituting $\PS(x)$ for $x$, we
find that the following inequality is true for all $x$:
\begin{align}
\nonumber
  g_{i(j)}((\PS(x))_{i(j)}) & - g_{i(j)}((x_{j+1})_{i(j)}) + \langle
  (\PS(x)-x_{j+1})_{i(j)}, \nabla_{i(j)} f(\hx_{j}) \rangle \\
& \geq
  \frac{\Lmax}{\gamma}\langle (\PS(x)-x_{j+1})_{i(j)},
  (\Delta_j)_{i(j)} \rangle.
\label{eqn_pre_1}
\end{align}
From the definition of $\Lmax$, and using the notation
\eqnok{eq:defdelta}, we have
\begin{align*}
f(x_{j+1}) \leq f(x_j) + \langle \nabla_{i(j)}f(x_j), -(\Delta_j)_{i(j)} \rangle + {\Lmax\over 2}|(\Delta_j)_{i(j)}|^2,
\end{align*}
or equivalently,
\begin{align}
\langle \nabla_{i(j)}f(x_j), (\Delta_j)_{i(j)} \rangle \leq f(x_j) - f(x_{j+1})  + {\Lmax\over 2}|(\Delta_j)_{i(j)}|^2.
\label{eqn_pre_2}
\end{align}
From 
the definition of $\bar{x}_{j+1}$ in \eqref{eqn_xj1}, we have
\[
\bar{x}_{j+1}=\arg\min_x \, \langle \nabla f(\hx_{j}), x- x_j\rangle +
    {\Lmax \over 2\gamma}\|x-x_j\|^2 + g(x),
\]
so, using (41) from \citep{Tseng08} again, we have
\begin{align}
g(x) - g(\bar{x}_{j+1}) + \left\langle x-\bar{x}_{j+1},  \nabla f(\hx_{j}) + \frac{\Lmax}{\gamma}(\bar{x}_{j+1}-x_j) \right\rangle \geq 0,\quad \forall \, x.
\label{eqn_pre_3}
\end{align}
We now define 
\begin{align}
\Delta_j:=x_j-\bar{x}_{j+1},
\label{eq:Delta}
\end{align} 
and note that this definition is consistent with $(\Delta_j)_{i(j)}$
defined in~\eqref{eq:defdelta}. From \eqref{xbarij}, we have
\begin{align}
\E_{i(j)} (\|x_{j+1}-x_j\|^2) = {1\over n} \|\bar{x}_{j+1}- x_j\|^2.
\label{eq:expdel}
\end{align}

Recalling that the indices in $K(j)$ are sorted in the increasing
order from smallest (oldest) iterate to largest (newest) iterate, we
use $K(j)_t$ to denote the $t$-th smallest entry in $K(j)$. 
For $T=0,1,\dotsc,|K(j)|$, we define
\[
\hx_{j, T} := \hx_{j} + \sum_{t=1}^T(x_{K(j)_t + 1} - x_{K(j)_t}).
\]
\vspace{10mm}
We have the following relations:
\begin{align*}
\hx_{j} &= \hx_{j, 0}\\
x_{j} &= \hx_{j, |K(j)|}\\
x_j-\hx_{j} &= \sum_{t=0}^{|K(j)|-1}(\hx_{j, t+1} - \hx_{j, t})\\
\nabla f(x_j)-\nabla f(\hx_{j})&=\sum_{t=0}^{|K(j)|-1}(\nabla f(\hx_{j, t+1}) - \nabla f(\hx_{j, t})).
\end{align*}
{\rc Furthermore, we have}
\begin{align}
\nonumber &\quad \|\nabla f(x_j) - \nabla f(\hx_{j})\| 
\\
& 
=\left\|\sum_{t=0}^{|K(j)|-1} (\nabla f(\hx_{j, t}) - \nabla f(\hx_{j, t+1}))\right\| \nonumber\\
\quad & 
\leq \sum_{t=0}^{|K(j)|-1}\|\nabla f(\hx_{j, t}) - \nabla f(\hx_{j, t+1})\| \nonumber\\
\quad &
\leq \Lres \sum_{t=0}^{|K(j)|-1}\|\hx_{j, t} - \hx_{j, t+1}\| \nonumber\\
\quad &
= \Lres \sum_{t=1}^{|K(j)|} \|x_{K(j)_t}-x_{K(j)_t+1}\| \nonumber \\
\quad & 
= \Lres \sum_{d\in K(j)} \|x_{d+1}-x_{d}\|, 
\label{eqn_gbound}
\end{align}
where the second inequality holds because $\hx_{j,t}$ and $\hx_{j,
  t+1}$ differ in only a single coordinate.
\medskip


\subsection{Proof of Theorem~\ref{thm_2}}

\begin{proof}
We prove \eqref{eqn_thm2_1} by induction. First, note that for any
vectors $a$ and $b$, we have
\begin{align*}
\|a\|^2 - \|b\|^2 = & 2\|a\|^2 -(\|a\|^2+\|b\|^2) 
\\ \leq & 
2\|a\|^2 -
2\langle a, b\rangle \\
= &
2\langle a, a-b \rangle  \\
\leq &
2\|a\|\|b-a\|.
\end{align*}
Thus for all $j$, we have
\begin{align}
\nonumber \|x_{j-1}-\bar{x}_j\|^2 & - \|x_j-\bar{x}_{j+1}\|^2 \\
\leq & 
2\|x_{j-1}-\bar{x}_j\| {\|x_j-\bar{x}_{j+1}-x_{j-1}+\bar{x}_j\|}.
\label{eqn_proof2_1}
\end{align}
The second factor in the r.h.s. of \eqref{eqn_proof2_1} is bounded as follows:
\begin{align}
\nonumber
&\|x_j-\bar{x}_{j+1}-x_{j-1}+\bar{x}_j\|\\
\nonumber
&= \Bigg\|x_j-\PO\left(x_j-{\gamma \over \Lmax} \nabla f(\hx_{j})\right)-\\
\nonumber & \quad 
\left(x_{j-1}-\PO\left(x_{j-1}-{\gamma \over \Lmax} \nabla f(\hx_{j-1})\right)\right)\Bigg\|\\
\nonumber
&\leq \|x_j- x_{j-1}\| + \\
\nonumber
&\quad
\left\|\PO\left(x_j-{\gamma \over \Lmax} \nabla f(\hx_{j})\right)-\PO\left(x_{j-1}-{\gamma \over \Lmax} \nabla f(\hx_{j-1})\right)\right\|\\
\nonumber
&\leq 2\|x_j-x_{j-1}\| + {\gamma \over \Lmax} \|\nabla f(\hx_{j}) - \nabla f(\hx_{j-1})\| 
\\
\nonumber
&\quad \left(\text{by the nonexpansive property of $\PO$}\right) \\
\nonumber
& = 2\|x_j-x_{j-1}\| + {\gamma \over \Lmax} \|\nabla f(\hx_{j}) - \nabla f(x_j) + \nabla f(x_j)  - \nabla f(x_{j-1})
\\ \nonumber
&\quad + \nabla f(x_{j-1}) - \nabla f(\hx_{j-1})\| \\
\nonumber
&\leq 2\|x_j-x_{j-1}\| + {\gamma \over \Lmax} \big(\|\nabla f(\hx_{j}) - \nabla f(x_j)\|
 + \|\nabla f(x_j)  - \nabla f(x_{j-1})\| 
\\ \nonumber & \quad\quad  
+ \|\nabla f(x_{j-1}) - \nabla f(\hx_{j-1})\| \big) \\
\nonumber
&\leq \left(2+{\Lambda\gamma}\right)\|x_j-x_{j-1}\| + {\gamma \over \Lmax} \|\nabla f(\hx_{j}) - \nabla f(x_j)\| 
\\ \nonumber
& \quad
+ {\gamma \over \Lmax}\|\nabla f(x_{j-1}) - \nabla f(\hx_{j-1})\| \\
\nonumber
&\leq \left(2+{\Lambda\gamma}\right)\|x_j-x_{j-1}\| + {\Lambda\gamma}\sum_{d\in K(j)}\|x_d - x_{d+1}\| 
\\ &\quad
+ {\Lambda\gamma}\sum_{d\in K(j-1)}\|x_d - x_{d+1}\| \quad (\text{from \eqref{eqn_gbound}})
\label{eqn_proof2_2} \\
\nonumber
&\leq \left(2+{\Lambda\gamma}\right)\|x_j-x_{j-1}\| + {\Lambda\gamma}\sum_{d=j-\tau}^{j-1}\|x_d - x_{d+1}\| + {\Lambda\gamma}\sum_{d=j-1-\tau}^{j-2}\|x_d - x_{d+1}\| \\
&\leq \left(2+2{\Lambda\gamma}\right)\|x_j-x_{j-1}\| + 2{\Lambda\gamma}\sum_{d=j-1-\tau}^{j-2}\|x_d - x_{d+1}\|,
\label{eqn_proof2_2add}
\end{align}
where 
the fourth inequality uses $\|\nabla f(x_j) - \nabla f(x_{j-1})\| \leq
\Lres \|x_j- x_{j-1}\|$, since $x_j$ and $x_{j-1}$ differ in just one
component.

We set $j=1$, and note that $K(0)=\emptyset$ and $K(1)\subset
\{0\}$. In this case, we obtain a bound from \eqnok{eqn_proof2_2}
\begin{align*}
&\|x_1-\bar{x}_{2}+x_{0}-\bar{x}_1\| \leq \left(2+{\Lambda\gamma}\right)\|x_1-x_{0}\| + {\Lambda\gamma} \|x_1-x_0\| = \left(2+{2\Lambda\gamma}\right)\|x_1-x_{0}\|.
\end{align*}
By substituting this bound in \eqref{eqn_proof2_1} and setting $j=1$,
and taking expectations, we obtain 
\begin{align}
\nonumber
\E(\|x_{0}-\bar{x}_1\|^2) - \E(\|x_1-\bar{x}_{2}\|^2) & \leq
2\E(\|x_{0}-\bar{x}_1\| {\|x_1-\bar{x}_{2}-x_{0}+\bar{x}_1\|}) \\
 \label{eq:crap13}
& \leq
\left(4+{4\Lambda\gamma}\right)\E(\|\bar{x}_1-x_{0}\|\|x_1-x_{0}\|).
\end{align}
For any positive scalars $\mu_1$, $\mu_2$, and $\alpha$, we have
\beq \label{eq:gtrick} \mu_1\mu_2 \leq {1\over 2}(\alpha \mu_1^2 +
\alpha^{-1} \mu_2^2).  
\eeq 
It follows that 
\begin{align}
\nonumber
&\E (\|x_{j}-x_{j-1}\| \|\bar{x}_j -x_{j-1}\|) \leq {1\over 2}\E(n^{1/2}\|x_{j}-x_{j-1}\|^2 + n^{-1/2}\|\bar{x}_j -x_{j-1}\|^2)\\
\nonumber
&\quad = {1\over 2}\E(n^{1/2}\E_{i(j-1)}(\|x_{j}-x_{j-1}\|^2) + n^{-1/2}\|\bar{x}_j -x_{j-1}\|^2)\\
\nonumber
&\quad = {1\over 2}\E\left(n^{-1/2}\|\bar{x}_{j}-x_{j-1}\|^2 + n^{-1/2}\|\bar{x}_j -x_{j-1}\|^2\right)\quad (\text{from \eqref{eq:expdel}})\\
\label{eqn_proof2_4}
&\quad = n^{-1/2}\E \|\bar{x}_j -x_{j-1}\|^2.
\end{align}
By taking $j=1$ in \eqref{eqn_proof2_4}, and substituting in
\eqnok{eq:crap13}, we obtain
\[
 \E(\|x_{0}-\bar{x}_1\|^2) - \E(\|x_1-\bar{x}_{2}\|^2) \leq n^{-1/2}\left(4+{4\Lambda\gamma}\right)\E \|\bar{x}_1 -x_{0}\|^2,
\]
which implies that
\[
 \E(\|x_{0}-\bar{x}_1\|^2) \leq \left(1- {4 + {4\gamma \Lambda}\over \sqrt{n}}\right)^{-1} \E(\|x_1-\bar{x}_{2}\|^2) \leq \rho  \E(\|x_1-\bar{x}_{2}\|^2).
\]
To see the last inequality, one only needs to verify that
\[
\rho^{-1} \leq 1-{4+ 4\gamma \Lambda \over \sqrt{n}} \;
\Leftrightarrow \;  \gamma \leq \frac{\sqrt{n}(1-\rho^{-1}) - 4}{4\Lambda},
\]
where the last inequality follows from the second bound for $\gamma$ in \eqref{eq:boundgammac}. We have thus shown that \eqref{eqn_thm2_1} holds for $j=1$.

To take the inductive step, we assume that
\eqref{eqn_thm2_1} holds up to index $j-1$. We have for $j-1-\tau \leq
d \leq j-2$ and any $\beta>0$ (using \eqref{eq:gtrick} again) that
\begin{alignat*}{2}
\nonumber
&\E (\|x_{d}-x_{d+1}\| \|\bar{x}_j -x_{j-1}\|) && \\
\nonumber
&\quad \leq {1\over 2}\E (n^{1/2}\beta\|x_{d}-x_{d+1}\|^2 + n^{-1/2}\beta^{-1}\|\bar{x}_j -x_{j-1}\|^2) && \\
\nonumber
&\quad = {1\over 2} \E (n^{1/2}\beta\E_{i(d)}(\|x_{d}-{x}_{d+1}\|^2) + n^{-1/2}\beta^{-1}\|\bar{x}_j -x_{j-1}\|^2) && \\
\nonumber
&\quad = {1\over 2} \E (n^{-1/2}\beta\|x_{d}-\bar{x}_{d+1}\|^2 + n^{-1/2}\beta^{-1}\|\bar{x}_j -x_{j-1}\|^2) \quad && (\text{from \eqref{eq:expdel}})\\
\nonumber
&\quad \leq {1\over 2} \E (n^{-1/2}\beta\rho^{j-1-d}\|x_{j-1}-\bar{x}_{j}\|^2 + n^{-1/2}\beta^{-1}\|\bar{x}_j -x_{j-1}\|^2) 
\\ \nonumber
&\quad \quad
(\text{by the inductive hypothesis}).
\end{alignat*}
Thus by setting $\beta=\rho^{(d+1-j)/2}$, we obtain
\begin{align}
\E (\|x_{d}-x_{d+1}\| \|\bar{x}_j -x_{j-1}\|)
 \leq {\rho^{(j-1-d)/2}\over n^{1/2}} \E \left(\|\bar{x}_j -x_{j-1}\|^2\right).
\label{eqn_proof2_3}
\end{align}

By substituting \eqref{eqn_proof2_2add} into \eqref{eqn_proof2_1} and
taking expectation on both sides of~\eqref{eqn_proof2_1}, we obtain
\begin{align*}
&\E( \|x_{j-1}-\bar{x}_j\|^2) - \E(\|x_j-\bar{x}_{j+1}\|^2)\\
\leq & 2\E(\|\bar{x}_j -x_{j-1}\| \|\bar{x}_j-\bar{x}_{j+1}+x_j-x_{j-1}\|)\\
\leq & 2\E\left(\|\bar{x}_j -x_{j-1}\| \left(\left(2+2\Lambda\gamma\right)\|x_j-x_{j-1}\| + 2\Lambda\gamma \sum_{d=j-1-\tau}^{j-2} \|x_{d} -x_{d+1}\|\right)\right)\\
= & \left(4+4\Lambda\gamma\right)\E(\|\bar{x}_j -x_{j-1}\| \|x_j-x_{j-1}\|) + 4\Lambda\gamma \sum_{d=j-1-\tau}^{j-2} \E (\|\bar{x}_j -x_{j-1}\|\|x_{d} -x_{d+1}\|) \\
\leq & n^{-1/2}(4+4\Lambda\gamma)\E(\|\bar{x}_j-x_{j-1}\|^2) 
\\ 
&\quad \quad
+ n^{-1/2}4\Lambda\gamma\E(\|x_{j-1}-\bar{x}_j\|^2)\sum_{d=j-1-\tau}^{j-2}\rho^{(j-1-d)/2}
\quad(\text{from \eqref{eqn_proof2_4} and \eqref{eqn_proof2_3}}) \\
\leq & n^{-1/2}(4+4\Lambda\gamma)\E(\|\bar{x}_j-x_{j-1}\|^2) + n^{-1/2}4\Lambda\gamma\E(\|x_{j-1}-\bar{x}_j\|^2)\sum_{t=1}^{\tau}\rho^{t/2}\\
= & n^{-1/2}\left(4+4\Lambda\gamma(1+\theta)\right) \E (\|{x}_{j-1}-\bar{x}_j\|^2),
\end{align*}
where the last equality follows from the definition of $\theta$ in \eqref{eq:defpsic}. It follows that
\begin{align*}
\E(\|x_{j-1}-\bar{x}_j\|^2) \leq & \left(1-n^{-1/2}\left(4+4\Lambda\gamma(1+\theta)\right)\right)^{-1} \E
(\|x_j-\bar{x}_{j+1}\|^2)  \\
\leq &
\rho \E (\|x_j-\bar{x}_{j+1}\|^2).
\end{align*}
To see the last inequality, one only needs to verify that
\[
\rho^{-1} \leq 1-{4+{4\gamma \Lambda(1+\theta)}\over \sqrt{n}} \;
\Leftrightarrow \;  \gamma \leq \frac{\sqrt{n}(1-\rho^{-1}) - 4}{4\Lambda(1+\theta)},
\]
and the last inequality is true because of the upper bound of $\gamma$
in \eqref{eq:boundgammac}.  We have thus proved \eqref{eqn_thm2_1}.

Next we will show the expectation of the objective $F$ is
monotonically decreasing. We have by using the definition
\eqnok{eq:defdelta} and \eqref{xbarij} that
\begin{align}
\nonumber
&\E_{i(j)}  F(x_{j+1}) = \E_{i(j)} \left[f(x_j - (\Delta_j)_{i(j)} e_{i(j)} ) + g(x_{j+1})\right]\\
\nonumber
&\leq \E_{i(j)}\Bigg[f(x_j) + \langle \nabla_{i(j)}f(x_j), (\bar{x}_{j+1}-x_j)_{i(j)} \rangle + \frac{\Lmax}{2}\|(x_{j+1}-x_j)_{i(j)}\|^2 \\
\nonumber
& \quad + g_{i(j)}((x_{j+1})_{i(j)})  + \sum_{l\neq i(j)} g_l((x_{j+1})_l)\Bigg] \\
\nonumber
& = \E_{i(j)}\Bigg[f(x_j) + \langle \nabla_{i(j)}f(x_j), (\bar{x}_{j+1}-x_j)_{i(j)} \rangle + \frac{\Lmax}{2}\|(x_{j+1}-x_j)_{i(j)}\|^2 \\
\nonumber 
& \quad+ g_{i(j)}((x_{j+1})_{i(j)})  + \sum_{l\neq i(j)} g_l((x_{j})_l)\Bigg] \\
\nonumber
& = f(x_j)+ {n-1\over n}g(x_j)
+n^{-1}\left( \langle \nabla f(x_j), \bar{x}_{j+1}-x_j \rangle + \frac{\Lmax}{2}\|\bar{x}_{j+1}-x_j\|^2+g(\bar{x}_{j+1})\right),
\end{align}
where we used $\E_{i(j)}\sum_{l\neq i(j)}g_l(x_j)_l={n-1\over
  n}g(x_j)$ in the last equality. By adding and subtracting a term
involving $\nabla f(\hat{x}_j)$, we obtain
\begin{align}
\nonumber
&\E_{i(j)} F(x_{j+1}) 
\\ \nonumber
& \leq 
F(x_j)+{1\over n}\left(\langle \nabla f(\hx_{j}), \bar{x}_{j+1}-x_j \rangle + \frac{\Lmax}{2}\|\bar{x}_{j+1}-x_j\|^2 + g(\bar{x}_{j+1}) - g(x_j) \right) \\
\nonumber
&\quad\quad+{1\over n} \langle \nabla f(x_j)- \nabla f(\hx_{j}), \bar{x}_{j+1}-x_j\rangle \\
\nonumber
& \leq F(x_j)+ {1\over n}\left(\frac{\Lmax}{2}\|\bar{x}_{j+1}-x_j\|^2  - {\Lmax\over \gamma}\|\bar{x}_{j+1}-x_j\|^2\right) \\
\nonumber
&\quad\quad+ {1\over n} \langle \nabla f(x_j)- \nabla f(\hx_{j}), \bar{x}_{j+1}-x_j\rangle \quad\quad (\text{from \eqref{eqn_pre_3} with $x=x_j$})
\\
&= F(x_j)-  \left({\frac{1}{\gamma}-\frac{1}{2}}\right)\frac{\Lmax}{n}\|\bar{x}_{j+1}-x_j\|^2 +  {1\over n}\langle \nabla f(x_j)- \nabla f(\hx_{j}), \bar{x}_{j+1}-x_j\rangle.
\label{eqn_proof2_5}
\end{align}
Consider the expectation of the last term on the right-hand side of
this expression. We have
\begin{alignat}{2}
\nonumber
\E \langle & \nabla f(x_j)- \nabla f(\hx_{j}), \bar{x}_{j+1}-x_j\rangle && \\
\nonumber
& \leq \E \left( \| \nabla f(x_j)- \nabla f(\hx_{j})\| \|\bar{x}_{j+1}-x_j \| \right) && \\
\nonumber
&\leq \Lres \E \left( \sum_{d\in K(j)}\|x_{d+1}-x_d\|\|\bar{x}_{j+1}-x_j\|\right) \quad (\text{from \eqref{eqn_gbound}}) \\
\nonumber
&\leq \Lres\sum_{d=j-\tau}^{j-1}{\rho^{(j-d)/2}\over n^{1/2}}\E (\|x_{j} - \bar{x}_{j+1}\|^2) \quad (\text{from \eqref{eqn_proof2_3}, and replacing $j$ by $j+1$})\\
&\leq n^{-1/2}\Lres \theta \E(\|x_j-\bar{x}_{j+1}\|^2) \quad (\text{from \eqref{eq:defpsic}}).
\label{eqn_proof2_5_5}
\end{alignat}
By taking expectations on both sides of \eqref{eqn_proof2_5} and
substituting \eqref{eqn_proof2_5_5}, we obtain
\[
\E F(x_{j+1}) \leq \E F(x_j) - {1\over n}\left(\left({\frac{1}{\gamma}-\frac{1}{2}}\right)\Lmax-{\Lres \theta\over n^{1/2}}\right)\E\|\bar{x}_{j+1}-x_j\|^2.
\]
To see $\left({\frac{1}{\gamma}-\frac{1}{2}}\right)\Lmax-{\Lres
  \theta\over n^{1/2}}\geq 0$ or equivalently
$\left({\frac{1}{\gamma}-\frac{1}{2}}\right)-{\Lambda \theta\over
  n^{1/2}}\geq 0$, we note from \eqref{eq:defpsic} and
\eqref{eq:boundgammac} that
\begin{align*}
\gamma^{-1} \geq \psi \geq {1\over 2} + \frac{\Lambda \theta}{\sqrt{n}}.
\end{align*}
Therefore, we have proved the monotonicity of the expectation of the
objectives, that is,
\begin{equation} \label{eq:monF}
\E F(x_{j+1}) \leq \E  F(x_{j}),\quad j=0,1,2,\dotsc.
\end{equation} 

Next we prove the sublinear convergence rate for the constrained
smooth convex case in~\eqref{eqn_thm2_2}. We have
\begin{alignat}{2}
\nonumber
&\| x_{j+1}- \PS(x_{j+1})\|^2 \leq \|x_{j+1} - \PS(x_j)\|^2 && \\
\nonumber
& = \|x_j-(\Delta_j)_{i(j)}e_{i(j)} - \PS(x_j)\|^2 && \\
\nonumber
& = \|x_j - \PS(x_j)\|^2 + |(\Delta_j)_{i(j)}|^2 - 2  \langle (x_j - \PS(x_j))_{i(j)}, (\Delta_j)_{i(j)}\rangle && \\
\nonumber
& = \|x_j - \PS(x_j)\|^2 - |(\Delta_j)_{i(j)}|^2 - 2  \langle (x_j - \PS(x_j))_{i(j)}-(\Delta_j)_{i(j)}, (\Delta_j)_{i(j)}\rangle && \\
\nonumber
& = \|x_j - \PS(x_j)\|^2 - |(\Delta_j)_{i(j)}|^2 + 2  \langle \PS(x_j)-x_{j+1})_{i(j)}, (\Delta_j)_{i(j)} \rangle 
\quad(\text{from \eqref{eq:defdelta}}) \\
\nonumber
& \leq \|x_j - \PS(x_j)\|^2 - |(\Delta_j)_{i(j)}|^2 + && \\
\nonumber 
&\quad  \frac{2\gamma}{\Lmax} \left[\langle (\PS(x_j)-x_{j+1})_{i(j)}, \nabla_{i(j)} f(\hx_{j})\rangle + g_{i(j)}((\PS(x_j))_{i(j)}) - g_{i(j)}((x_{j+1})_{i(j)}) \right] 
\\ \nonumber
&\quad \quad
\quad (\text{from \eqref{eqn_pre_1}})
\\
\nonumber
& = \|x_j - \PS(x_j)\|^2 - |(\Delta_j)_{i(j)}|^2 + && \\
\nonumber
&\quad \frac{2\gamma}{\Lmax} \left[\langle (\PS(x_j)-x_j)_{i(j)}, \nabla_{i(j)} f(\hx_{j})\rangle + g_{i(j)}((\PS(x_j))_{i(j)}) - g_{i(j)}((x_{j+1})_{i(j)})  \right] + && \\
\nonumber
&\quad \quad \frac{2\gamma}{\Lmax}\left(\langle (\Delta_j)_{i(j)}, \nabla_{i(j)} f(x_{j})\rangle +
\langle (\Delta_j)_{i(j)}, \nabla_{i(j)} f(\hx_{j}) - \nabla_{i(j)} f(x_{j})\rangle \right) && \\
\nonumber
& \leq \|x_j - \PS(x_j)\|^2 - |(\Delta_j)_{i(j)}|^2 + && \\
\nonumber
&\quad  \frac{2\gamma}{\Lmax} \left[\langle(\PS(x_j)-x_j)_{i(j)}, \nabla_{i(j)} f(\hx_{j})\rangle + g_{i(j)}((\PS(x_j))_{i(j)}) - g_{i(j)}((x_{j+1})_{i(j)})  \right] + && \\
\nonumber
&\quad  \frac{2\gamma}{\Lmax}\bigg(f(x_j) -  f(x_{j+1})  + {\Lmax\over 2}|(\Delta_j)_{i(j)}|^2 +  \langle (\Delta_j)_{i(j)}, \nabla_{i(j)} f(\hx_{j}) - \nabla_{i(j)} f(x_{j})\rangle \bigg) 
\\ \nonumber
& \quad \quad \quad 
(\text{from \eqref{eqn_pre_2}})\\
\nonumber
& = \|x_j - \PS(x_j)\|^2 - (1-\gamma)|(\Delta_j)_{i(j)}|^2 + \frac{2\gamma}{\Lmax} \underbrace{\langle (\PS(x_j)-x_j)_{i(j)}, \nabla_{i(j)} f(\hx_{j})\rangle}_{T_1} +&& \\
\nonumber 
&\quad \frac{2\gamma}{\Lmax}\underbrace{\langle (\Delta_j)_{i(j)}, \nabla_{i(j)} f(\hx_{j}) - \nabla_{i(j)} f(x_{j}) \rangle}_{T_2} + &&
\\ \label{eqn_proof2_6} &\quad
\frac{2\gamma}{\Lmax}\underbrace{\left[f(x_j) -  f(x_{j+1}) + g_{i(j)}((\PS(x_j))_{i(j)}) - g_{i(j)}((x_{j+1})_{i(j)})  \right]}_{T_3}. &&
\end{alignat}
We now seek upper bounds on the quantities $T_1$, $T_2$, and $T_3$ in
the expectation sense. For simplicity, we construct a vector
$\bold{b}\in\R^{|K(j)|}$ with $\bold{b}_t = \|\hx_{j, t-1}-\hx_{j,
  t}\|$. We have from elementary arguments that
\begin{alignat}{2}
\nonumber
\E( & \|\bold{b}\|^2) = \sum_{t=0}^{|K(j)|-1} \E (\|\hx_{j, t}-\hx_{j, t+1}\|^2) = \sum_{t=1}^{|K(j)|} \E (\|x_{K(j)_t}-x_{K(j)_t+1}\|^2) && \\ 
\nonumber 
& =
\sum_{d\in K(j)} \E (\|x_{d}-x_{d+1}\|^2) = \frac{1}{n} \sum_{d\in K(j)} \E\|x_d-\bar{x}_{d+1}\|^2 
\leq \frac{1}{n}\sum_{d=j-\tau}^{j-1} \E\|x_d-\bar{x}_{d+1}\|^2 && \\ 
\nonumber 
& \leq \frac{1}{n} \sum_{t=1}^{\tau} \rho^{t}\E\|x_j-\bar{x}_{j+1}\|^2 \quad  (\text{from \eqref{eqn_thm2_1}})
\\ 
\label{eq:b}
& \leq {\theta' \over n} \E\|x_j-\bar{x}_{j+1}\|^2 \quad (\text{from \eqref{eq:defpsic}}).
\end{alignat}
For the expectation of $T_1$, defined in \eqref{eqn_proof2_6}, we have
\begin{align}
\nonumber 
\E( & T_1) = \E \left( (\PS(x_j)-x_j)_{i(j)} \nabla_{i(j)} f(\hx_{j}) \right) \\
\nonumber
& = n^{-1}\E \langle \PS(x_j)-x_j, \nabla f(\hx_{j})  \rangle \\
\nonumber
& = n^{-1}\E \langle \PS(x_j)-\hx_{j}, \nabla f(\hx_{j})  \rangle
 + n^{-1}\E \sum_{t=0}^{|K(j)|-1}\langle \hx_{j, t}-\hx_{j, t+1}, \nabla f(\hx_{j})  \rangle \\
\nonumber
& = n^{-1}\E \langle \PS(x_j)-\hx_{j}, \nabla f(\hx_{j})  \rangle  \\
\nonumber
& \quad + n^{-1}\E \sum_{t=0}^{|K(j)|-1} \left(\langle \hx_{j, t}-\hx_{j, t+1}, \nabla f(\hx_{j, t})  \rangle+ \langle \hx_{j, t}-\hx_{j, t+1}, \nabla f(\hx_{j}) - \nabla f(\hx_{j, t})  \rangle \right)
\\ \nonumber
& \leq n^{-1}\E (f^*_j- f(\hx_{j})) 
\\ \nonumber
& \quad
+ n^{-1} \E \sum_{t=0}^{|K(j)|-1}\left( f(\hx_{j, t})-f(\hx_{j, t+1}) + {\Lmax\over 2}\|\hx_{j, t}-\hx_{j, t+1}\|^2 \right) \\
\nonumber
& \quad +n^{-1} \E \sum_{t=0}^{|K(j)|-1} \langle \hx_{j, t}-\hx_{j, t+1}, \nabla f(\hx_{j}) - \nabla f(\hx_{j, t})  \rangle \quad (\text{from \eqref{eq:col}})
\\ \nonumber
& = n^{-1}\E (f^*_j- f(x_{j})) + {\Lmax\over 2n} \E \|\bold{b}\|^2 
\\ \nonumber
&\quad
+n^{-1} \E \sum_{t=0}^{|K(j)|-1} \langle \hx_{j, t}-\hx_{j, t+1}, \nabla f(\hx_{j}) - \nabla f(\hx_{j, t})  \rangle\\
\nonumber
& = n^{-1}\E (f^*_j- f(x_{j})) + {\Lmax\over 2n} \E \|\bold{b}\|^2 
\\ \nonumber
& \quad 
+n^{-1} \E \sum_{t=0}^{|K(j)|-1} \left\langle \hx_{j, t}-\hx_{j, t+1}, \sum_{t'=0}^{t-1}\nabla f(\hx_{j,t'}) - \nabla f(\hx_{j,t'+1})  \right\rangle \\
\nonumber
& \leq n^{-1}\E (f^*_j- f(x_{j})) + {\Lmax\over 2n} \E \|\bold{b}\|^2 
\\ \nonumber
& \quad 
+n^{-1} \E \sum_{t=0}^{|K(j)|-1}  \Lmax\left(\|\hx_{j, t}-\hx_{j, t+1}\|\sum_{t'=0}^{t-1} \|\hx_{j,t'}-\hx_{j,t'+1}\|\right) \\
\nonumber
& = n^{-1}\E (f^*_j- f(x_{j})) + {\Lmax\over 2n} \E \|\bold{b}\|^2 +n^{-1} \Lmax\E \sum_{t=0}^{|K(j)|-1} \left( \bold{b}_{t+1} \sum_{t'=0}^{t-1} \bold{b}_{t'+1} \right) \\
\nonumber
& = n^{-1}\E (f^*_j- f(x_{j})) + {\Lmax\over 2n} \E \|\bold{b}\|^2 + {\Lmax \over 2n} \E(\| \bold{b} \|_1^2 - \|\bold{b}\|^2) \\
\nonumber
& = n^{-1}\E (f^*_j- f(x_{j})) + {\Lmax \over 2n} \E( \| \bold{b}\|_1^2 ) \\
\nonumber
& \leq n^{-1}\E (f^*_j- f(x_{j})) + {\Lmax \tau \over 2n} \E(\|\bold{b}\|^2) \quad (\text{since $\|\bold{b}\|_1 \leq \sqrt{|K(j)|}\|\bold{b}\| \leq \sqrt{\tau}\|\bold{b}\|$})
\\
&\leq n^{-1}\E (f^*_j- f(x_{j})) + {\Lmax \tau\theta' \over 2n^2} \E(\|x_j-\bar{x}_{j+1}\|^2) \quad (\text{from \eqref{eq:b}}).
\label{eq:bT1}
\end{align}
For the expectation of $T_2$, we have \\
\vspace{-10mm}
\begin{alignat}{2}
\nonumber
&\E(T_2) =\E  (\Delta_j)_{i(j)} \left( \nabla_{i(j)} f(\hx_{j}) - \nabla_{i(j)} f(x_{j}) \right) && \\
\nonumber
&\quad= n^{-1}\E \langle \Delta_j, \nabla f(\hx_{j}) - \nabla f(x_j) \rangle && \\
\nonumber
&\quad\leq n^{-1} \E (\|\Delta_j\| \|\nabla f(\hx_{j}) - \nabla f(x_j)\|) && \\
\nonumber
&\quad\leq {\Lres \over n} \E \left( \sum_{d=j-\tau}^{j-1} \|\Delta_j\| \|x_{d} - x_{d+1}\| \right) \quad && (\text{from \eqref{eqn_gbound}})\\
\nonumber
&\quad= \frac{\Lres}{n} \E \left( \sum_{d=j-\tau}^{j-1} \|x_j-\bar{x}_{j+1}\| \|x_{d} - x_{d+1}\| \right) && \\
\nonumber
&\quad\leq \frac{\Lres}{n^{3/2}} \sum_{d=j-\tau}^{j-1} \rho^{(j-d)/2}\E\|x_j-\bar{x}_{j+1}\|^2 
\quad && (\text{from \eqref{eqn_proof2_3} with $j$ replacing $j-1$})\\
\label{eq:bT2}
&\quad \leq \frac{\Lres\theta}{n^{3/2}} \E\|x_j-\bar{x}_{j+1}\|^2 \quad 
&& (\text{from \eqref{eq:defpsic}}).
\end{alignat}
For $T_3$, let us look the expectation of several individual terms first
\[\E_{i(j)} g_{i(j)}((\PS(x_j))_{i(j)}) = n^{-1}g(\PS(x_j)) = n^{-1}g^*_j,\] and
\begin{align*}
\E_{i(j)} g_{i(j)}((x_{j+1})_{i(j)}) &=  \E_{i(j)} (g(x_{j+1})-g(x_j) + g_{i(j)}((x_j)_{i(j)})) 
\\ &=  
\E_{i(j)} g(x_{j+1}) - g(x_j) + n^{-1} g(x_j) 
\\ &= 
\E_{i(j)} g(x_{j+1}) - {n-1\over n}g(x_j).
\end{align*}
Now we take the expectation on $T_3$ and use the equalities above to obtain:
\begin{align}
\nonumber 
\E(T_3) &= \E f(x_j) - \E f(x_{j+1}) + \E g_{i(j)}((\PS(x_j))_{i(j)}) - \E g_{i(j)}((x_{j+1})_{i(j)}) 
\\
\label{eq:bT3}
& = \E f(x_j) - \E f(x_{j+1}) + n^{-1} \E g^*_j - \E g(x_{j+1}) + {n-1\over n} \E g(x_j). 
\end{align}
By substituting the upper bounds from \eqnok{eq:bT1}, \eqnok{eq:bT2}, and \eqnok{eq:bT3} into \eqref{eqn_proof2_6}, we obtain
\begin{align}
\nonumber
\E\|x_{j+1} & - \PS(x_{j+1})\|^2 \leq \E \|x_j - \PS(x_j)\|^2 - (1-\gamma) \E |(\Delta_j)_{i(j)}|^2 
\\ \nonumber & +
{2\gamma \over \Lmax}\left({1\over n}\E (f^*_j - f(x_j)) + {\Lmax \tau\theta' \over 2n^2} \E(\|x_j-\bar{x}_{j+1}\|^2) \right) 
\\ \nonumber & +
{2\gamma \over \Lmax}\left(\frac{\Lres\theta}{n^{3/2}} \E\|x_j-\bar{x}_{j+1}\|^2\right)
\\ \nonumber & +
{2\gamma \over \Lmax}\left(\E f(x_j) - \E f(x_{j+1}) + n^{-1} \E g^*_j - \E g(x_{j+1}) + {n-1\over n} \E g(x_j)\right).
\end{align}
By using
\[\E_{i(j)}(|(\Delta_j)_{i(j)}|^2) = n^{-1} \| x_j-\bar{x}_{j+1} \|^2,\]
it follows that
\begin{align}
\nonumber
\E & \|x_{j+1} - \PS(x_{j+1})\|^2 \leq \E \|x_j - \PS(x_j)\|^2 \\
\nonumber
& \quad - {1\over n}\left(1-\gamma - {\tau\theta' \over n}\gamma - {2\Lambda\theta \over n^{1/2}}\gamma\right)\E\|x_j-\bar{x}_{j+1}\|^2 \\
\nonumber
&\quad +\frac{2\gamma}{\Lmax n}(\E f^*_j- \E f(x_j) + \E g^*_j) 
\\ \nonumber
& \quad  
+ \frac{2\gamma}{\Lmax}(\E f(x_j) + {n-1\over n}\E g(x_j) -  \E f(x_{j+1}) - \E g(x_{j+1}))  \\
\nonumber
&\leq \E \|x_j - \PS(x_j)\|^2 + \frac{2\gamma}{\Lmax n}(\E f^*_j- \E f(x_j) + \E g^*_j) \\
\nonumber
&\quad 
+\frac{2\gamma}{\Lmax}\left(\E f(x_j) + {n-1\over n}\E g(x_j) -  \E f(x_{j+1}) - \E g(x_{j+1})\right) \\
\label{eqn_proof2_7_5}
&\leq \E \|x_j - \PS(x_j)\|^2 + \frac{2\gamma}{\Lmax n}(F^*- \E F(x_j)) + \frac{2\gamma}{\Lmax}(\E F(x_j) -  \E F(x_{j+1}) ).
\end{align}
In the second inequality, we were able to drop the term involving $\E
\|x_j-\bar{x}_{j+1}\|^2$ by using the fact that
\[
1-\gamma\left(1+\frac{\tau\theta'}{n} + \frac{\Lambda \theta}{\sqrt{n}}\right) = 1 - \gamma\psi \geq 0,
\]
which follows from the definition \eqnok{eq:defpsic} of $\psi$ and
from the first upper bound on $\gamma$ in \eqref{eq:boundgammac}.  It
follows from \eqref{eqn_proof2_7_5} that
\begin{align}
\nonumber
&\E\|x_{j+1} - \PS(x_{j+1})\|^2 + \frac{2\gamma}{\Lmax}(\E F(x_{j+1}) -F^*) \\
\label{eq:crap6}
&\quad \leq  \E \|x_j - \PS(x_j)\|^2 + \frac{2\gamma}{\Lmax}(\E F(x_j) - F^*) - \frac{2\gamma}{L_{\max}n}( \E F(x_j)- F^*). 
\end{align}
Defining
\begin{equation} \label{eq:defsj}
S_j:=\E (\|x_j-\PS(x_j)\|^2) + {2\gamma \over \Lmax}\E (F(x_j)-F^*),
\end{equation}
we have from \eqref{eq:crap6} that
\begin{equation} \label{eq:recursj}
S_{j+1}\leq S_j-{2\gamma\over {\Lmax}n}\E (F(x_j)-F^*),
\end{equation}
so by induction, we have
\begin{equation} \label{eq:sjbound}
S_{j+1} \leq S_0 - {2\gamma\over \Lmax n} \sum_{t=0}^{j}(\E F(x_t) - F^*) \le
 S_0 - {2\gamma (j+1) \over \Lmax n} (F(x_0)-F^*),
\end{equation}
where the second inequality follows from monotonicity of $\E F(x_j)$
\eqnok{eq:monF}. Note that
\[
S_0:=\|x_0-\PS(x_0)\|^2 + {2\gamma \over \Lmax} (F(x_0)-F^*).
\]
By substituting the definition of $S_{j+1}$ into \eqref{eq:sjbound},
we obtain
\begin{align*}
\E\|x_{j+1} & - \PS(x_{j+1})\|^2 + \frac{2\gamma}{\Lmax}(\E
    F(x_{j+1}) -F^*) + \frac{2\gamma (j+1)}{L_{\max}n}(\E F(x_{j+1})-
    F^*)\\ 
& \leq \| x_0-\PS(x_0)\|^2 +  \frac{2\gamma}{\Lmax}(F(x_0) - F^*).
\end{align*}
The sublinear convergence expression \eqnok{eqn_thm2_2} follows when
we drop the (nonnegative) first term on the left-hand side of this
expression, and rearrange.

Finally, we prove the linear convergence rate \eqnok{eqn_thm2_3} for
the optimally strongly convex case. All bounds proven above continue
to hold, and we make use the optimal strong convexity property in
\eqref{eq:esc}:
\begin{align*}
F(x_j) - F^* \geq \frac{l}{2}\|x_j - \PS(x_j)\|^2.
\end{align*}
By using this result together with some elementary manipulation, we
obtain
\begin{align}
\nonumber
F(x_j) - F^* &= \left(1- \frac{\Lmax}{l\gamma+\Lmax}\right)(F(x_j) - F^*) +   \frac{\Lmax}{l\gamma+\Lmax}(F(x_j) - F^*)\\
\nonumber
&\geq \left(1- \frac{\Lmax}{l\gamma+\Lmax}\right)(F(x_j) - F^*) +   \frac{L_{\max}l}{2(l\gamma+\Lmax)}\|x_j - \PS(x_j)\|^2\\
\label{eqn_proof2_9}
&= \frac{L_{\max}l}{2(l\gamma+\Lmax)}\left(\|x_j - \PS(x_j)\|^2 + \frac{2\gamma}{\Lmax}(F(x_j) - F^*)\right).
\end{align}
By taking expectations of both sides in this expression, and comparing with
\eqref{eq:defsj}, we obtain
\[
\E (F(x_j) - F^*) \ge \frac{\Lmax l}{2(l \gamma+\Lmax)} S_j.
\]
By substituting into \eqnok{eq:recursj}, we obtain
\begin{align*}
S_{j+1}  \le & S_j - \left(\frac{2 \gamma}{\Lmax n}\right) \frac{\Lmax l}{2(l \gamma+\Lmax)} S_j  
\\=& 
\left( 1-\frac{l \gamma}{n(l\gamma + \Lmax)} \right) S_j 
\\ \le & 
\left( 1-\frac{l \gamma}{n(l\gamma + \Lmax)} \right)^{j+1} S_0,
\end{align*}
where the last inequality follows from induction over $j$. We obtain
\eqref{eqn_thm2_3} by substituting the definition \eqref{eq:defsj} of
$S_j$.
\end{proof}

\subsection{Proof of Corollary~\ref{co:thm_2}}

\begin{proof}
Note that for $\rho$ defined by~\eqref{eq:choicerhoc}, and using
\eqnok{eq:boundtauc}, we have
\begin{align} 
\nonumber
\rho^{(1+\tau)/2}= & \left(1+ {4e\Lambda (\tau + 1) \over
  \sqrt{n}}\right)^{1+\tau} =  \left(\left(1+ {4e\Lambda (\tau + 1) \over
  \sqrt{n}}\right)^{\sqrt{n}\over 4e\Lambda(\tau +1)}\right)^{{4e\Lambda (\tau+1)^2 \over \sqrt{n}}} 
\\ &
\leq e^{{4e\Lambda(\tau+1)^2
    \over \sqrt{n}}} \leq e.
\label{eq:rhot1}
\end{align}
Thus from the definition of $\psi$~\eqref{eq:defpsic}, we have that
\begin{alignat*}{2}
\psi & = 1+\frac{\tau\theta'}{n} + \frac{2\Lambda\theta}{\sqrt{n}} 
\\ \nonumber & \leq 
1+\frac{\tau^2\rho^{\tau}}{n} + \frac{2\Lambda\tau\rho^{\tau/2}}{\sqrt{n}} \quad 
&& \left(\text{from $\theta = \sum_{t=1}^{\tau}\rho^{t/2} \leq \tau \rho^{\tau/2}$ and $\theta' = \sum_{t=1}^{\tau}\rho^{t} \leq \tau \rho^{\tau}$}\right)
\\ & \leq
1+\frac{\tau^2e^2}{n} + \frac{2\Lambda\tau e}{\sqrt{n}} && \left(\text{from \eqref{eq:rhot1}}\right)
\\
&\leq
1+{1\over 16} + {1\over 2} \leq 2,
\end{alignat*}
where for the second-last inequality we used \eqref{eq:boundtauc} to obtain
\[
\frac{\Lambda \tau e}{\sqrt{n}} \le \frac{\Lambda \tau e}{4e \Lambda (\tau+1)^2} \le \frac14, \quad
\frac{\tau^2 e^2}{n} = \left( \frac{\tau e}{\sqrt{n}} \right)^2 \le
\left( \frac{\Lambda \tau e}{\sqrt{n}} \right) \le \frac{1}{16}.
\]
Thus, the steplength parameter choice $\gamma=1/2$ satisfies the first
bound in~\eqref{eq:boundgammac}. To show that the second bound in
\eqnok{eq:boundgammac} holds also, we have
\begin{alignat*}{2}
&\frac{\sqrt{n}(1-\rho^{-1})-4}{4(1+\theta){\Lambda}}  &&
\\ &
\geq \frac{\sqrt{n}(1-\rho^{-1})}{4(1+\theta){\Lambda}} - {1\over 2}  
&& \text{(from $\theta \geq 1$ and $\Lambda \geq 1$)}
\\ 
&
\geq \frac{\sqrt{n}(1-\rho^{-1/2})}{4(1+\theta){\Lambda}} - {1\over 2}  
&&
\\ 
&
= \frac{\sqrt{n}(\rho^{1/2}-1)}{4(1+\theta)\rho^{1/2}{\Lambda}} - {1\over 2}
\\ &
\geq \frac{\sqrt{n}(\rho^{1/2}-1)}{4(\tau+1)\rho^{(\tau+1)/2}{\Lambda}} - {1\over 2}
\quad && \left(\text{from $(1+\theta)\rho^{1/2}\leq (1+\tau\rho^{\tau/2})\rho^{1/2} \leq (1+\tau)\rho^{(\tau+1)/2}$}\right)
\\ &
\geq \frac{4e\Lambda(\tau+1)}{4e(\tau+1){\Lambda}} - {1\over 2} 
&& \text{(from \eqref{eq:choicerhoc} and \eqref{eq:rhot1})}
\\ &
\geq 1 - {1\over 2} ={1\over 2}. &&
\end{alignat*}
We can thus set $\gamma=1/2$, and by substituting this choice
into~\eqref{eqn_thm2_3}, we obtain~\eqref{eqn_thm_2_good_c}. We
obtain~\eqref{eqn_thm_3_good_c} by making the same substitution
into~\eqref{eqn_thm2_2}.
\end{proof}

{
\bibliographystyle{siam}
\bibliography{reference}

\begin{thebibliography}{10}

\bibitem{AgarwalD12}
{\sc A.~Agarwal and J.~C. Duchi}, {\em Distributed delayed stochastic
  optimization}, in Proceedings of the Conference on Decision and Control,
  2012, pp.~5451--5452.

\bibitem{Ani00a}
{\sc M.~Anitescu}, {\em Degenerate nonlinear programming with a quadratic
  growth condition}, SIAM Journal on Optimization, 10 (2000), pp.~1116--1135.

\bibitem{Avron13arXiv}
{\sc H.~{Avron}, A.~{Druinsky}, and A.~{Gupta}}, {\em Revisiting asynchronous
  linear solvers: Provable convergence rate through randomization}, in
  Proceedings of the {IEEE} International Parallel and Distributed Processing
  Symposium, May 2014.

\bibitem{BeckT09}
{\sc A.~Beck and M.~Teboulle}, {\em A fast iterative shrinkage-thresholding
  algorithm for linear inverse problems}, SIAM Journal on Imaging Sciences, 2
  (2009), pp.~183--202.

\bibitem{Beck13}
{\sc A.~Beck and L.~Tetruashvili}, {\em On the convergence of block coordinate
  descent type methods}, SIAM Journal on Optimization, 23 (2013),
  pp.~2037--2060.

\bibitem{Bertsekas89}
{\sc D.~P. Bertsekas and J.~N. Tsitsiklis}, {\em Parallel and Distributed
  Computation: Numerical Methods}, Prentice Hall, 1989.

\bibitem{Boyd11}
{\sc S.~Boyd, N.~Parikh, E.~Chu, B.~Peleato, and J.~Eckstein}, {\em Distributed
  optimization and statistical learning via the alternating direction method of
  multipliers}, Foundations and Trends in Machine Learning, 3 (2011),
  pp.~1--122.

\bibitem{Bradley11}
{\sc J.~K. Bradley, A.~Kyrola, D.~Bickson, and C.~Guestrin}, {\em Parallel
  coordinate descent for {L}1-regularized loss minimization}, in International
  Conference on Machine Learning, 2011.

\bibitem{CortesVapnik95}
{\sc C.~Cortes and V.~Vapnik}, {\em Support vector networks}, Machine Learning,
   (1995), pp.~273--297.

\bibitem{Cotter11}
{\sc A.~Cotter, O.~Shamir, N.~Srebro, and K.~Sridharan}, {\em Better mini-batch
  algorithms via accelerated gradient methods}, in Advances in Neural
  Information Processing Systems, vol.~24, 2011, pp.~1647--1655.

\bibitem{Duchi12}
{\sc J.~C. Duchi, A.~Agarwal, and M.~J. Wainwright}, {\em Dual averaging for
  distributed optimization: Convergence analysis and network scaling}, IEEE
  Transactions on Automatic Control, 57 (2012), pp.~592--606.

\bibitem{Elsner92}
{\sc L.~Elsner, I.~Koltracht, and M.~Neumann}, {\em Convergence of sequential
  and asynchronous paracontractions nonlinear paracontractuions}, Numerische
  Mathematik, 62 (1992), pp.~305--316.

\bibitem{FacSS13}
{\sc F.~Facchinei, S.~Sagratella, and G.~Scutari}, {\em Flexible parallel
  algorithms for big data optimization}, technical report, Department of
  Computer, Control, and Management Engineering, University of Rome "La
  Sapienza", November 2013.
\newblock arXiv:1311.2444v1.

\bibitem{Fercoq13b}
{\sc O.~Fercoq and P.~Richt{\'a}rik}, {\em Accelerated, parallel, and proximal
  coordinate descent}, technical report, School of Mathematics, University of
  Edinburgh, 2013.
\newblock arXiv: 1312.5799.

\bibitem{Fercoq13a}
\leavevmode\vrule height 2pt depth -1.6pt width 23pt, {\em Smooth minimization
  of nonsmooth functions by parallel coordinate descent}, technical report,
  School of Mathematics, University of Edinburgh, 2013.
\newblock arXiv:1309.5885.

\bibitem{Ferris94}
{\sc M.~C. Ferris and O.~L. Mangasarian}, {\em Parallel variable distribution},
  SIAM Journal on Optimization, 4 (1994), pp.~815--832.

\bibitem{Frommer00}
{\sc A.~Frommer and D.~B. Szyld}, {\em On asynchronous iterations}, Journal of
  Computational and Applied Mathematics, 123 (2000), pp.~201--216.

\bibitem{Ma12}
{\sc D.~Goldfarb and S.~Ma}, {\em Fast multiple-splitting algorithms for convex
  optimization}, SIAM Journal on Optimization, 22 (2012), pp.~533--556.

\bibitem{Hoffman52}
{\sc A.~J. Hoffman}, {\em On approximate solutions of systems of linear
  inequalities}, Journal of Research of the National Bureau of Standards, 49
  (1952), pp.~263--265.

\bibitem{LaiYin13}
{\sc M.~Lai and W.~Yin}, {\em Augmented {L}1 and nuclear-norm models with a
  globally linearly convergent algorithm}, SIAM Journal on Imaging Sciences, 6
  (2013), pp.~1059--1091.

\bibitem{LiuWright13}
{\sc J.~{Liu}, S.~J. {Wright}, C.~{R{\'e}}, V.~{Bittorf}, and S.~{Sridhar}},
  {\em An asynchronous parallel stochastic coordinate descent algorithm},
  technical report, Computer Sciences Department, University of
  Wisconsin-Madison, February 2014.
\newblock arXiv: 1311.1873.

\bibitem{Liu14arXivAsyRK}
{\sc J.~{Liu}, S.~J. {Wright}, and S.~{Sridhar}}, {\em An asynchronous parallel
  randomized {Kaczmarz} algorithm}, technical report, Computer Sciences
  Department, University of Wisconsin-Madison, 2014.
\newblock arXiv: 1401.4780.

\bibitem{LuXiao13}
{\sc Z.~Lu and L.~Xiao}, {\em On the complexity analysis of randomized
  block-coordinate descent methods}, Technical Report MSR-TR-2013-53, Microsoft
  Research, May 2013.
\newblock arXiv:1305.4723.

\bibitem{LuoTseng92}
{\sc Z.-Q. Luo and P.~Tseng}, {\em On the convergence of the coordinate descent
  method for convex differentiable minimization}, Journal of Optimization
  Theory and Applications, 72 (1992), pp.~7--35.

\bibitem{Mangasarian95}
{\sc O.~L. Mangasarian}, {\em Parallel gradient distribution in unconstrained
  optimization}, SIAM Journal on Optimization, 33 (1995), pp.~916--1925.

\bibitem{Necoara13b}
{\sc I.~{Necoara} and D.~{Clipici}}, {\em Efficient parallel coordinate descent
  algorithm for convex optimization problems with separable constraints:
  application to distributed {MPC}}, technical report, Automation and Systems
  Engineering Department, University Politechnica Bucharest, 2013.
\newblock arXiv: 1302.3092.

\bibitem{Necoara13a}
{\sc I.~{Necoara} and A.~{Patrascu}}, {\em A random coordinate descent
  algorithm for optimization problems with composite objective function and
  linear coupled constraints}, technical report, Automation and Systems
  Engineering Department, University Politechnica Bucharest, 2013.
\newblock arXiv: 1302.3074.

\bibitem{Nemirovski09}
{\sc A.~Nemirovski, A.~Juditsky, G.~Lan, and A.~Shapiro}, {\em Robust
  stochastic approximation approach to stochastic programming}, SIAM Journal on
  Optimization, 19 (2009), pp.~1574--1609.

\bibitem{nesterov2004introductory}
{\sc Y.~Nesterov}, {\em Introductory Lectures on Convex Optimization: A Basic
  Course}, Kluwer Academic Publishers, 2004.

\bibitem{Nesterov12}
\leavevmode\vrule height 2pt depth -1.6pt width 23pt, {\em Efficiency of
  coordinate descent methods on huge-scale optimization problems}, SIAM Journal
  on Optimization, 22 (2012), pp.~341--362.

\bibitem{Hogwild11nips}
{\sc F.~Niu, B.~Recht, C.~R{\'e}, and S.~J. Wright}, {\em Hogwild: A lock-free
  approach to parallelizing stochastic gradient descent}, Advances in Neural
  Information Processing Systems, 24 (2011), pp.~693--701.

\bibitem{Yin13}
{\sc Z.~Peng, M.~Yan, and W.~Yin}, {\em Parallel and distributed sparse
  optimization}, tech. report, Department of Mathematics, UCLA, 2013.

\bibitem{Richtarik11}
{\sc P.~{Richt{\'a}rik} and M.~{Tak{\'a}{\v c}}}, {\em Iteration complexity of
  randomized block-coordinate descent methods for minimizing a composite
  function}, Mathematical Programing, Series A,  (2012).
\newblock (Published Online).

\bibitem{Richtarik12arXiv}
\leavevmode\vrule height 2pt depth -1.6pt width 23pt, {\em Parallel coordinate
  descent methods for big data optimization}, technical report, Mathematics
  Department, University of Edinburgh, 2012.
\newblock arXiv: 1212.0873.

\bibitem{Saha13}
{\sc A.~Saha and A.~Tewari}, {\em On the nonasymptotic convergence of cyclic
  coordinate descent methods}, SIAM Journal on Optimization, 23 (2013),
  pp.~576--601.

\bibitem{ScherrerTHH12}
{\sc C.~Scherrer, A.~Tewari, M.~Halappanavar, and D.~Haglin}, {\em Feature
  clustering for accelerating parallel coordinate descent}, in Advances in
  Neural Information Processing, vol.~25, 2012, pp.~28--36.

\bibitem{Shalev-Shwartz2013}
{\sc S.~Shalev-Shwartz and T.~Zhang}, {\em Accelerated mini-batch stochastic
  dual coordinate ascent}, in Advances in Neural Information Processing
  Systems, C.~J.~C. Burges, L.~Bottou, M.~Welling, Z.~Ghahramani, and K.~Q.
  Weinberger, eds., vol.~26, 2013, pp.~378--385.

\bibitem{Shamir2013icml}
{\sc O.~Shamir and T.~Zhang}, {\em Stochastic gradient descent for non-smooth
  optimization: Convergence results and optimal averaging schemes}, in
  Proceedings of the International Conference on Machine Learning, 2013.

\bibitem{Sridhar2013nips}
{\sc S.~{Sridhar}, V.~{Bittorf}, J.~{Liu}, C.~{Zhang}, C.~{R{\'e}}, and S.~J.
  {Wright}}, {\em An approximate efficient solver for {LP} rounding}, in
  Advances in Neural Information Processing Systems, vol.~26, 2013.

\bibitem{Tseng01}
{\sc P.~Tseng}, {\em Convergence of a block coordinate descent method for
  nondifferentiable minimization}, Journal of Optimization Theory and
  Applications, 109 (2001), pp.~475--494.

\bibitem{Tseng08}
{\sc P.~Tseng}, {\em On accelerated proximal gradient methods for
  convex-concave optimization}, technical report, University of Washington,
  2008.

\bibitem{TseY06}
{\sc P.~Tseng and S.~Yun}, {\em A coordinate gradient descent method for
  nonsmooth separable minimization}, Mathematical Programming, Series {B}, 117
  (2009), pp.~387--423.

\bibitem{TseY07a}
\leavevmode\vrule height 2pt depth -1.6pt width 23pt, {\em A coordinate
  gradient descent method for linearly constrained smooth optimization and
  support vector machines training}, Computational Optimization and
  Applications, 47 (2010), pp.~179--206.

\bibitem{WangLin13}
{\sc P.-W. Wang and C.-J. Lin}, {\em Iteration complexity of feasible descent
  methods for convex optimization}, technical report, Department of Computer
  Science, National Taiwan University, 2013.

\bibitem{WriNF08a}
{\sc S.~J. Wright, R.~D. Nowak, and M.~A.~T. Figueiredo}, {\em Sparse
  reconstruction by separable approximation}, IEEE Transactions on Signal
  Processing, 57 (2009), pp.~2479--2493.

\end{thebibliography}
}

\end{document}